\documentclass[10pt,a4paper]{amsart}

\usepackage[marginpar=2cm,ignoremp,margin=3cm]{geometry}

\RequirePackage{doi}
\usepackage{hyperref}

\usepackage{xfrac}
\newcommand{\half}{{\sfrac{1\!}{2}}}
\newcommand{\quart}{{\sfrac{1\!}{4}}}

\allowdisplaybreaks

\usepackage{amssymb}

\def\equationautorefname~#1\null{Equation~(#1)\null}

\usepackage{mathtools}
\usepackage{thmtools}

\usepackage{todonotes}

\declaretheorem[
style=plain,
name=Theorem,
refname={Theorem,Theorems},
Refname={Theorem,Theorems}
]{Thm}
\declaretheorem[
style=plain,
name=Proposition,
numberlike=Thm,
refname={Proposition,Propositions},
Refname={Proposition,Propositions}
]{Prop}
\declaretheorem[
style=plain,
name=Lemma,
numberlike=Thm,
refname={Lemma,Lemmas},
Refname={Lemma,Lemmas}
]{Lem}
\declaretheorem[
style=plain,
name=Corollary,
numberlike=Thm,
refname={Corollary,Corollaries},
Refname={Corollary,Corollaries}
]{Cor}

\DeclareMathOperator{\Li}{Li}
\renewcommand{\Im}{\operatorname{Im}}
\renewcommand{\Re}{\operatorname{Re}}

\renewcommand{\epsilon}{\varepsilon}
\newcommand{\ii}{\mkern 1mu \mathrm{i} \mkern 1mu}
\let\eps\varepsilon
\newcommand{\od}{\mathrm{od}}
\newcommand{\ev}{\mathrm{ev}}

\def\+{\!+\!}
\def\-{\!-\!}

\makeatletter
\def\paragraph{\@startsection{paragraph}{4}%
	{0em}{0.75ex \@plus 0.5ex \@minus 0.25ex}{-\fontdimen2\font}%
	{%
		\bfseries
}}
\makeatother

\let\overlineO\overline
\renewcommand{\overline}[1]{\overlineO{\mathclap{\phantom{I}}#1}}

\usepackage{stackengine}
\newcommand\brabar{\scalebox{.3}{(\,}\raisebox{-2.1pt}{--}\scalebox{.3}{\,)}}
\newcommand\brabarB{\scalebox{.3}{(\,}\raisebox{-1.5pt}{--}\scalebox{.3}{\,)}} 
\newcommand\bratw{{\scalebox{.4}{\normalfont (\,}\raisebox{-6pt}{%
$\widetilde{\smash{\phantom{I_{a,b}}}}$%
}\scalebox{.4}{\normalfont \,)}}} 

\makeatletter
\@namedef{subjclassname@2020}{%
	\textup{2020} Mathematics Subject Classification}
\makeatother

\begin{document}
	
	\title[{On evaluations of $S(\stackon[.1pt]{$2$}{\brabarB}, 1, \ldots, 1, \stackon[.1pt]{$1$}{\brabarB})$ and $T(\stackon[.1pt]{$2$}{\brabarB}, 1, \ldots, 1, \stackon[.1pt]{$1$}{\brabarB})$}]{On the evaluations of multiple \( S \) and \( T \) values \\ of the form 
	\( S(\stackon[.1pt]{$2$}{\brabar}, 1, \ldots, 1, \stackon[.1pt]{$1$}{\brabar}) \) and \( T(\stackon[.1pt]{$2$}{\brabar}, 1, \ldots, 1, \stackon[.1pt]{$1$}{\brabar}) \) 
	\\[0.5ex] {\footnotesize \normalfont Answers to questions of Xu, Yan, and Zhao}}
	\author{Steven Charlton}
	\date{7 March 2024}
	
	\hypersetup{
			pdftitle={On the evaluation of multiple S and T values of the form S(±2, 1, …, ±1) and T(±2, 1, …, ±1)},
			pdfauthor = {Steven Charlton},
	}
	
	\address{Max Planck Institute for Mathematics, Vivatsgasse 7, Bonn 53111, Germany}
	\email{charlton@mpim-bonn.mpg.de}
	
	\keywords{(alternating) multiple $T$ values, multiple $S$ values, multiple $t$ values, multiple zeta values, weighted sums, iterated integrals, multiple polylogarithms, generalised doubling identity, shuffle antipode, regularisation, special values, generating series, arctangent integrals}
	\subjclass[2020]{11M32}
	
	\begin{abstract}
		Xu, Yan and Zhao showed that in even weight, the multiple \( T \) value \( T(2, 1, \ldots, 1, \overline{1}) \) is a polynomial in \( \log(2), \pi \), Riemann zeta values, and Dirichlet beta values.  Based on low-weight examples, they conjectured that \( \log(2) \) does not appear in the evaluation.  We show that their conjecture is correct, and in fact follows largely from various earlier results of theirs.  More precisely, we derive explicit formulae for \( T(2, 1, \ldots, 1, \overline{1}) \) in even weight and \( S(2, 1, \ldots, 1, \overline{1}) \) in odd weight via generating series calculations.  We also resolve another conjecture of theirs on the evaluations of \( T(\overline{2}, 1, \ldots, 1, \overline{1}) \), \( S(\overline{2}, 1, \ldots, 1, 1) \), and \( S(\overline{2}, 1, \ldots, 1, \overline{1}) \) in even weight, by way of calculations involving Goncharov's theory of iterated integrals and multiple polylogarithms.	
	\end{abstract}
	
	\maketitle
	
	\section{Introduction and statement}
	
	The alternating multiple mixed values (AMMV's \cite{xu2022alternating}, but appearing already in \cite{xu2024multipleT}) are an extension and generalisation of many related objects (multiple $t$ values \cite{hoffman2019oddvariant}, multiple $T$ values \cite{kaneko2020level2}, multiple $S$ values \cite{xu2022variants}, and alternating versions thereof, as well as the multiple zeta values and their alternating versions \cite{hoffman1992harmonic,zagier1994zeta,broadhurst1996enumeration}), putting them all into a more uniform framework.  Each of these objects is defined as a nested sum over \( m_1 > m_2 > \cdots > m_r \), whose summation indices satisfy some particular parity conditions.  The alternating multiple mixed value allows the these parities to be specified arbitrarily.\medskip
	
	For any tuple of parities \( (\eps_i) \in \{ \pm 1 \}^r \), and a tuple of signs \( (\sigma_i) \in \{ \pm 1 \}^r \), the alternating multiple mixed value is defined in \cite{xu2024multipleT,xu2022alternating} by
	\begin{equation}\label{eqn:MM:def}
	\begin{aligned}
		& M^{\eps_1,\ldots,\eps_r}_{\sigma_1,\ldots,\sigma_r}(s_1,\ldots,s_r) \\
		& \coloneqq \sum_{m_1 > \cdots > m_r} \frac{(1 + \eps_1(-1)^{m_1}) \sigma_1^{(2m_1 + 1 - \eps_1)/4} \cdots (1 + \eps_r(-1)^{m_r}) \sigma_r^{(2m_r + 1 - \eps_r)/4}}{m_1^{s_1} \cdots m_r^{s_r}}
	\end{aligned}
	\end{equation}
	In particular, when \( \eps_i = 1 \), the combination \( (1 + \eps_i(-1)^{n_i}) \) requires \( n_i \) to be even, and when \( \eps_i = -1 \), it requires \( n_i \) to be odd, for a non-zero contribution to the sum.  By abuse of notation, one can write \( \eps_i = 1 \) as \( \eps_i = \ev \), and \( \eps_i = -1 \) as \( \eps_i = \od \) to emphasise the parity restriction these choices force on the summation indices. (Note: the summation and sign conventions differ in much of the literature \cite{kaneko2020level2,xu2021duality,xu2022variants,xu2024multipleT}, we adapt all formulae to the current $>$-convention, established in \cite{xu2022alternating}.)\medskip
	
	For simplicity of exposition, we shall take by definition, the alternating multiple \( T \)-values and \( S \)-values to be
	\begin{align*}
		T^{\sigma_1,\ldots,\sigma_r}(s_1,\ldots,s_r) \coloneqq M_{\sigma_1,\ldots,\sigma_r}^{\ldots, \ev, \od, \ev, \od}(s_1,\ldots,s_r) \\[1ex]
		S^{\sigma_1,\ldots,\sigma_r}(s_1,\ldots,s_r) \coloneqq M_{\sigma_1,\ldots,\sigma_r}^{\ldots, \od, \ev, \od, \ev}(s_1,\ldots,s_r) \,.
	\end{align*}
	In each case the parities flip between \( \od \) and \( \ev \), but multiple \( T \) values have the final parity (corresponding to the last/smallest summation index \( m_r \)) being odd, whilst multiple \( S \) values have it being even.  In contrast, the alternating multiple \( t \) values have all parities odd, giving (possibly up to a sign, and power of 2, depending on conventions)
	\[
		t^{\sigma_1,\ldots,\sigma_r}(s_1,\ldots,s_r) = \frac{1}{2^r} M^{\od, \ldots, \od}_{\sigma_1,\ldots,\sigma_r}(s_1,\ldots,s_r) \,,
	\]
	Likewise, the alternating multiple zeta values are obtained by taking all parities even (possibly up to a sign, and power of 2)
	\begin{equation}\label{eqn:zeta}
		\zeta^{\sigma_1,\ldots,\sigma_r}(s_1,\ldots,s_r) = \frac{2^{s_1 + \cdots + s_r}}{2^r} M^{\ev, \ldots, \ev}_{\sigma_1,\ldots,\sigma_r}(s_1,\ldots,s_r) \,.
	\end{equation}
	As is commonly done, we will indicate the signs in multiple \( T \) and \( S \) values (this also applies to multiple zeta and \( t \) values elsewhere in the literature) by writing a bar over the argument \( s_i \) if and only if \( \sigma_i = -1 \), so for example
	\[
		T(4, \overline{3}, \overline{2},1) = T^{1,-1,-1,1}(4, 3, 2, 1) \,. 
	\]
	Likewise, we will write \( \{a\}_n \) as shorthand for the string \( a, \ldots, a \), with exactly \( n \) repetitions of a.
	\medskip
	
	In  \cite{xu2022alternating}, Xu, Yan, and Zhao investigated algebraic and structural properties of the AMMV's, in particular establishing shuffle, and stuffle relations, regularisation behaviours, integral representations, and duality and parity results.  The authors also investigated some special values, and their relations to certain arctangent integrals.  In particular, they considered 
	\begin{align*}
		S(2, \{1\}_{2m-2}, \overline{1}) &= 2^{2m} \!\!\! \sum_{\substack{n_1 > n_2 > \cdots > n_{2m} \\ \text{$n_1,n_3,\ldots,n_{2m\-1}$ odd} \\ \text{$n_2,n_4,\ldots,n_{2m}$ even}}} \!\!\! \frac{(-1)^{n_{2m}/2}}{n_1^2 n_2 \cdots n_{2m}}
		\\
		T(2, \{1\}_{2m-1}, \overline{1}) &= - 2^{2m+1} \!\!\! \sum_{\substack{n_1 > n_2 > \cdots > n_{2m\+1} \\ \text{$n_1,n_3,\ldots,n_{2m\+1}$ odd} \\ \text{$n_2,n_4,\ldots,n_{2m}$ even}}} \!\!\! \frac{(-1)^{(n_{2m+1}-1)/2}}{n_1^2 n_2 \cdots n_{2m\+1}} \,.
	\end{align*}
	In Corollary 4.8 \cite{xu2022alternating}, they conclude \vspace{1ex}
	\[
		S(2, \{1\}_{2m-2}, \overline{1}) \,,\, T(2, \{1\}_{2m-1}, \overline{1}) \in \mathbb{Q}[\log(2), \pi, \zeta(2), \beta(2), \zeta(3), \beta(3), \ldots] \,,\vspace{1ex}
	\]
	where here
	\[
		\zeta(s) \coloneqq \sum_{n=1}^\infty \frac{1}{n^s} \,, \quad \Re(s) > 1 \,, \quad \beta(s) \coloneqq \sum_{n=1}^\infty \frac{(-1)^{n-1}}{(2n-1)^s} \,, \quad \Re(s) > 0 \,,
	\]
	are the Riemann zeta function, and the Dirichlet beta function respectively.  Based on the second and fourth examples from this list,
	\begin{align*}
		S(2, \overline{1}) &= \frac{7}{2} \zeta(3) - \pi G - \frac{\pi^2}{4} \log(2) \,, \\
		T(2, 1, \overline{1}) &= - 6 \beta(4) + 3 \zeta(2) G \,, \\
		S(2, 1, 1, \overline{1}) &= \frac{31}{5} \zeta(5) - \frac{15}{8} \zeta(4) \log(2) - \frac{63}{32} \zeta(2)\zeta(3) - \pi \beta(4) \,, \\
		T(2, 1,  1, 1, \overline{1}) &= \frac{15}{4} \zeta(4) G + 3 \zeta(2) \beta(4) - 10 \beta(6) \,,
	\end{align*}
	with \( G = \beta(2) \), the Catalan constant, they pose the following question.  	(These expressions can be obtained from the database of level \( 4 \) coloured MZV's tabulated by Au \cite{au2020evaluation}, or verified numerically using techniques or programs for numerical evaluation of multiple zeta values such as \texttt{polylogmult} in \texttt{pari/gp} \cite{PARI2} or the routine \verb|zeta({k1,...,kr},{s1,...,sr})| in the \texttt{GiNaC} \cite{GINAC} interactive shell.)%
	\medskip

	\noindent {\bf Question {\normalfont (Question 1, \cite[p. 18]{xu2022alternating})}.}  Is it true that\vspace{1ex}
	\[
		 T(2, \{1\}_{2m-1}, \overline{1}) \in \mathbb{Q}[\pi, \zeta(2), \beta(2), \zeta(3), \beta(3), \ldots] \,,\vspace{1ex}
	\]
	so that in particular \( \log(2) \) does not appear in the evaluation? \medskip
	
	We show an affirmative answer to this question, by giving an explicit generating series identity for both the multiple \( T \) value evaluation, and the multiple \( S \) value evaluation.  Unexpectedly, but pleasantly, the proof of this follows (more or less) directly from the proof of the authors result in \cite[Corollary 4.8]{xu2022alternating}, and earlier work of some of these authors \cite{xu2021duality,xu2022variants,xu2024multipleT}. \medskip
	
	Introduce the generating series of the relevant multiple \( S \) values and multiple \( T \) values,
	\begin{align*}
	 	E(z) \coloneqq \sum_{m=0}^\infty S(2, \{1\}_{2m}, \overline{1}) z^{2m+2} \,, \quad
	 	F(z) \coloneqq \sum_{m=0}^\infty T(2, \{1\}_{2m+1}, \overline{1}) z^{2m+3} \,.
	\end{align*}
	Define the following generating series of odd zeta values (as in \cite{zagier2012evaluation}), and even beta values
	\[
		A(z) \coloneqq \sum_{k=1}^\infty \zeta(2k+1) z^{2k} \,, \quad D(z) \coloneqq \sum_{k=1}^\infty \beta(2k) z^{2k-1} \,.
	\]
	The first main result of this note is then as follows.	
	\begin{Thm}\label{thm:qn1:eval}
		The following generating series identity holds
		\begin{align*}
			E(z) + \ii F(z) =  
			\begin{aligned}[t]
			& - \Big( A\Big( \frac{z}{4} \Big) - 3 A \Big( \frac{z}{2} \Big) + 2 A(z) + \log(2) \Big) \cdot \frac{\pi z}{2}  \tan\Big( \frac{\pi z}{2} \Big)  \\
			& - z \Big(  \frac{1}{2} A'\Big( \frac{z}{2} \Big) - 2 A'(z) - \pi D(z) \Big)  \\
			& + \ii  \Big(  2 z \beta(2) - 2 z D'(z) + 2 \cdot \frac{\pi z}{2} \tan\Big(\frac{\pi z}{2} \Big)  D(z) \Big) \,.
			\end{aligned}
		\end{align*}
	\end{Thm}
	One immediately sees that \( \log(2) \) does not appear in the evaluation of \( T(2, \{1\}_{2m+1}, \overline{1}) \) which is encapsulated in the imaginary part.  By using that
	\[
		\sum_{n=1}^\infty 2 (1 - 2^{-2n}) \zeta(2n) z^{2n} = \frac{\pi z}{2} \tan\Big( \frac{\pi z}{2} \Big) \,,
	\]
	and extracting the real and imaginary parts (equivalently the odd and even powers of \( z \)), we obtain the following more precise evaluations.
	\begin{Cor}
		The following evaluations hold
		\begin{align*}
			S(2, \{1\}_{2m-2}, \overline{1}) & = \begin{aligned}[t]
					&  - m \big( 2^{-2m+1} - 4 \big) \zeta(2m+1) - 2 (1 - 2^{-2m}) \log(2)  \zeta(2m)  \\
					& -\pi \beta(2m) \,\, + \!\! \sum_{\substack{a+b = m \\ a, b \geq 1}} \!\! 2(1 - 2^{-2a})\big( 4^{-2b}  - 3 \cdot 2^{-2b} + 2 \big) \zeta(2a) \zeta(2b+1) \,, \end{aligned} \\
		\shortintertext{and}
			T(2, \{1\}_{2m-1}, \overline{1}) &= -2(2m+1)\beta(2m+2) + \!\! \sum_{\substack{a + b = m+1 \\ a, b \geq 1}} \!\! 4 (1 - 2^{-2a}) \zeta(2a) \beta(2b) \,. 
		\end{align*}
	\end{Cor}
	These formulae agree with the examples provided earlier, they answer Question~1 from \cite[p. 18]{xu2022alternating} affirmatively, and  refine the characterisation given in \cite[Corollary 4.8]{xu2022alternating}.  Specifically
	\begin{align*}
		S(2, \{1\}_{2m-2}, \overline{1}) & \in \mathbb{Q} \langle \pi^{2m - 2k} \zeta(2k+1) \mid k = 1, \ldots, m \rangle \oplus \mathbb{Q}\langle \pi^{2m} \log(2) , \pi \beta(2m) \rangle  \,, \\
		T(2, \{1\}_{2m-1}, \overline{1}) & \in \mathbb{Q}\langle \pi^{2m+2-2k} \beta(2k) \mid k = 2, \ldots, 2m+2 \rangle \,.
	\end{align*}\medskip
	
	The authors also considered the further special values (obtained by changing the sign of the first, and/or last argument)
		\begin{align*}
	S(\overline{2}, \{1\}_{2m-1}, \overline{1}) &= 2^{2m+1} \!\!\! \sum_{\substack{n_1 > n_2 > \cdots > n_{2m\+1} \\ \text{$n_1,n_3,\ldots,n_{2m\+1}$ even} \\ \text{$n_2,n_4,\ldots,n_{2m}$ odd}}} \!\!\! \frac{(-1)^{n_1/2} (-1)^{n_{2m+1}/2}}{n_1^2 n_2 \cdots n_{2m\+1}} \,,
	\\
	S(\overline{2}, \{1\}_{2m}) &= 2^{2m+1} \!\!\! \sum_{\substack{n_1 > n_2 > \cdots > n_{2m\+1} \\ \text{$n_1,n_3,\ldots,n_{2m\+1}$ even} \\ \text{$n_2,n_4,\ldots,n_{2m}$ odd}}} \!\!\! \frac{(-1)^{n_{1}/2}}{n_1^2 n_2 \cdots n_{2m\+1}} \,,
	\\
	T(\overline{2}, \{1\}_{2m-1}, \overline{1}) &= 2^{2m+1} \!\!\! \sum_{\substack{n_1 > n_2 > \cdots > n_{2m\+1} \\ \text{$n_1,n_3,\ldots,n_{2m\+1}$ odd} \\ \text{$n_2,n_4,\ldots,n_{2m}$ even}}} \!\!\! \frac{(-1)^{(n_1-1)/2}(-1)^{(n_{2m+1}-1)/2}}{n_1^2 n_2 \cdots n_{2m\+1}} \,.
	\end{align*}
	Based on more low-weight evaluations, they posed another Question, as follows.\medskip
	
	\noindent {\bf Question {\normalfont (Question 2, \cite[p.~21]{xu2022alternating})}.}  Is it true that
	\[
		T(\overline{2}, \{1\}_{2m-1}, \overline{1}) \in \mathbb{Q}[\pi, \zeta(2), \beta(2), \zeta(3), \beta(3), \ldots] \,?
	\]
	Furthermore, is it true that \( S(\overline{2}, \{1\}_{2m}), S(\overline{2}, \{1\}_{2m-1}, \overline{1}) \) can be expressed as a polynomial in \( \log(2), \pi, \zeta(2k+1), \beta(2p) \) and double zetas of the form \( \zeta(\overline{2\ell+1}, 1)\), \( \ell, p, k \in \mathbb{Z}_{>0} \)? \medskip
	
	We show that the answer is again affirmative for the first two parts, but that the evaluation for \( S(\overline{2}, \{1\}_{2m-1}, \overline{1}) \) begins to invoke more complicated alternating double zetas values from weight 8 onwards, and so the answer is \emph{negative}.  More precisely we establish the following formula for the MTV, and the following characterisations of the two MSV's.
	
	\begin{Thm}\label{thm:qn2:mtv}
		For any \(m \geq 1 \), the following evaluation holds
		\[
			T(\overline{2}, \{1\}_{2m-1}, \overline{1}) = 
			 \sum_{\substack{r+p+2k=2m \\ r,p,k\geq0}} \!\! 4 (-1)^{r+k} \beta(r+1) \beta(p+1) \cdot \frac{1}{(2k)!} \Big( \frac{\pi}{2} \Big)^{2k} \,.
		\]
		So in particular, 
		\(
			T(\overline{2}, \{1\}_{2m-1}, \overline{1}) \in \mathbb{Q}[\pi^2, \beta(\mathrm{even})] 
		\).
	\end{Thm}

	\begin{Thm}\label{thm:qn2:msv}
		For any \( m \geq 1 \), we have
		\begin{align*}
			S(\overline{2}, \{1\}_{2m}) &\in \mathbb{Q}[\log(2), \pi, \zeta(\mathrm{odd}), \beta(\mathrm{even})] \oplus \mathbb{Q} \cdot \zeta(\overline{2m+1}, 1) \,, \\
			S(\overline{2}, \{1\}_{2m-1}, \overline{1}) &\in \mathbb{Q}[\log(2), \pi, \zeta(\mathrm{odd}), \beta(\mathrm{even})] \oplus \mathbb{Q} \cdot \zeta(\overline{2m+1}, 1) \oplus \mathbb{Q} \cdot W_{m} \,,
		\end{align*}
		where \( W_m \) denotes the following weighted sum of double alternating MZV's
		\[
			W_m \coloneqq \sum_{\substack{p + q = 2m+2 \\ p, q \geq 1}} \frac{1}{2^p} \zeta(\overline{\mkern 1mu p\mkern 1mu }, q) \,, \quad \text{ with } \quad  \zeta(\overline{\mkern 1mu a \mkern 1mu}, b) \coloneqq \sum_{n > m > 0} \frac{(-1)^a}{n^a m^b} \,.
		\]
	\end{Thm}
	In showing the second part of this theorem, we first obtain an evaluation for \( S(\overline{2}, \{1\}_{2m-1}, \overline{1}) \) in terms of \( \zeta(\overline{2m+1}, 1) \), Riemann zeta values, Dirichlet beta values and the alternating multiple \( t \) value \( t(\overline{1}, \overline{2m+1}) \).  It is the latter object which generates the more complicated alternating MZV contribution \( W_m \); in \autoref{prop:tm1modd:eval} we give an explicit evaluation for the double $t$ value, obtained via a generalised doubling relation \cite[\S14.2.5]{zhaoFunctions}, \cite[\S4]{mzvDM}.  (Note: A generalised doubling relation was already used in \cite{CK2242} to reduce \( \zeta(\overline{\ev},\overline{\ev}) \) to non-alternating double zeta values.  It would be interesting to investigate whether \( t(\overline{\od},\overline{\od}) \) can always be reduced to alternating double zeta values, or non-alternating double $t$ values, in a similar way.)
	In \autoref{cor:qn2:msvpart1} and \autoref{cor:qn2:msvpart2}, we give explicit formulae for the evaluations of these MSV's.
	\medskip
	
	The paper is structured as follows.  In \autoref{sec:thm1} we prove \autoref{thm:qn1:eval}, by making explicit the results established by Xu, Yan and Zhao leading up to Corollary 4.8 in \cite{xu2022alternating}.  Then we will recall some details about Goncharov's setup of iterated integrals in \autoref{sec:goncharov}, before establishing \autoref{thm:qn2:mtv} in \autoref{sec:thm2a}, and \autoref{thm:qn2:msv} in \autoref{sec:thm2b} and \autoref{sec:thm2c} by computations with iterated integrals, and properties of multiple polylogarithms.
	
	\paragraph{\normalfont\bf Acknowledgements}  The results of \autoref{thm:qn1:eval} were already directly communicated to the authors of \cite{xu2022alternating}, who did not expect Question 1 could so readily be resolved from their earlier results.  At their suggestion, I prepared a stand-alone note to formally present the solution, which they could then reference.  During the preparation of the note I was also able to resolve their Question~2, via different techniques, leading to the current---more substantial---text.
	
	I am grateful to the Max Planck Institute for Mathematics, Bonn, for support, hospitality and excellent working conditions during the preparation of this paper.
	
	\section{Proof of \autoref{thm:qn1:eval}}\label{sec:thm1}
	
	We follow the proof of the Lemmas, Propositions and Theorems leading up to Corollary 4.8 in \cite{xu2022alternating}, making all of the stages more explicit, and providing generating series expressions where necessary.	The key idea of the proof of Corollary 4.8  \cite{xu2022alternating} is to give two different expressions for the following ``arctangent over $x$'' integral:
	\[
	\int_0^1 \frac{\arctan^{m}(x)}{x} \, \mathrm{d}x \,.
	\]
	The first expression goes via \( T(\overline{2}, \{1\}_{p-1}) \) which is known to evaluate \cite[Theorems 3.16, 3.17]{xu2024multipleT} in terms of zeta values and beta values.  The second expression involves the values \( T(2, \{1\}_{2p-1}, \overline{1}) \) and \( S(2, \{1\}_{2p-2}, \overline{1}) \) which appear in the statement of the Corollary, and the values \( T(2, \{1\}_{2p}) \) and \( S(2, \{1\}_{2p+1}) \) which are already understood and evaluate \cite[Equation (3.17)]{xu2022variants} in terms of zeta values and beta values. Corollary 4.8 \cite{xu2022alternating} then follows from equating these two expressions. \medskip
	
	For the first expression: from the definition of arctan, and properties of iterated integrals, we have
	\begin{align*}
		\arctan^r(x) &= r! \int_{x >  t_1 > \cdots > t_r > 0} \frac{\mathrm{d}t_1}{1 + t_1^2}  \cdots \frac{\mathrm{d}t_r}{1 + t_r^2}   \\
		& = r! \sum_{n_1 > n_2 > \cdots > n_r > 0} \frac{(-1)^{n_1 - r} x^{2n_1 - r}}{(2n_1 - r)(2n_2 - r+1) \cdots (2n_r-1)} \,.
	\end{align*}
	By multiplying by \( \frac{1}{x} \), integrating this term-wise over \( 0 < x < 1 \), and using the definition of multiple \( T \) values, Xu, Yan, and Zhao then establish in Proposition 4.7 \cite{xu2022alternating} the equality
	\begin{equation}\label{eqn:arctanoverx:eval}
		\int_0^1 \frac{\arctan^r(x)}{x} \, \mathrm{d}x = (-1)^{\lfloor (r+1)/2 \rfloor} \frac{r!}{2^r} T(\overline{2}, \{1\}_{r-1}) \,.
	\end{equation}
	
	So form the generating series (note the shifted factorial in the first series),
	\[
	Q(z) \coloneqq \sum_{r=1}^\infty \frac{z^r}{(r-1)!} \int_0^1
 \frac{\arctan^r(x)}{x} \, \mathrm{d}x \,,
  \quad
  	G_{\bar2\{1\}}^T(z) \coloneqq \sum_{m=0}^\infty T(\overline{2}, \{1\}_m) z^{m+1}  \,.
	\]
	(Later, we will need other generating series of particular families of multiple \( T \) values and multiple \( S \) values; the notation above is chosen to help reminder the reader, and the author, which generating series is which at a glance.)  Using the result that if
	\(
		f(y) = \sum_n a_n y^n 
	\),
	then
	\[
		\sum_n (-1)^{\lfloor (r+1)/2 \rfloor} n a_n y^{n-1} = -\frac{1+\ii}{2} f'(-\ii y) - \frac{1-\ii}{2} f'(\ii y) \,,
	\]
	we deduce from \autoref{eqn:arctanoverx:eval} the following lemma.
	\begin{Lem}\label{lem:Q:expr1}
		The following generating series expression holds
		\[
			Q(z) = -\frac{z}{2} \Big( \frac{1+\ii}{2} {G_{\bar2\{1\}}^T}'\Big({-}\frac{\ii z}{2}\Big) + \frac{1-\ii}{2} {G_{\bar2\{1\}}^T}'\Big(\frac{\ii z}{2}\Big) \Big)  \,.
		\]
	\end{Lem}
	The first task then is to give an explicit generating series expression for \( T(\overline{2}, \{1\}_{r-1}) \). 
	
	\paragraph{Evaluation of \( T(\overline{2}, \{1\}_{r-1}) \), with weighted sum formulae}  By combining results from \cite{xu2021duality,xu2022variants,xu2024multipleT}, we can derive such a generating series expression. \medskip
	
	From \cite[Proposition 3.21]{xu2024multipleT}, we can deal with the even weight case, via the evaluation
	\begin{equation}\label{eqn:tm111m1}
		T(\overline{2}, \{1\}_{2p-2}) = -\sum_{1 \leq k \leq j \leq p} 4 \beta(2k) \cdot \frac{(1 - 2^{1 - (2j-2k)}) \zeta(2j-2k)}{2^{2j-2k}} \cdot  \frac{(-1)^{p-j}\pi^{2p-2j}}{(2p - 2j + 1)!} \,.
	\end{equation}
	The expression given in \cite{xu2024multipleT} involves alternating MtV's, and alternating MZV's, which we have directly expressed via beta values and zeta values for simplicity.  Likewise we have already substituted the definition \( \alpha_n \coloneqq \frac{\pi^{2n}}{(2n+1)!} \), c.f. Theorem 1.1 \cite{xu2024multipleT}.  (We have also adapted it to the current sign and summation order conventions.)
	
	Taking \( \sum_{p=1}^\infty z^{2p-1} \bullet \) of both sides of \autoref{eqn:tm111m1} gives the following lemma directly.
	\begin{Lem}\label{lem:gs:Tm2111}
		The following generating series identity holds
		\[
			\frac{1}{2} \big( G_{\bar2\{1\}}^T(z) - G_{\bar2\{1\}}^T(-z) \big) = \sum_{p=1}^\infty T(\overline{2}, \{1\}_{2p-2}) z^{2p-1} = - 2 \cos\Big( \frac{\pi z}{2} \Big) D(z) \,.
		\]
	\end{Lem}

	Now recall the weighted-sums from \cite{xu2021duality,xu2022variants,xu2024multipleT}
	\[
		W(k,r) \coloneqq \sum_{i_1 + \cdots + i_r = k} T(\overline{i_1}, i_2, \ldots, i_r) \,.
	\]
	From Theorem 2.3 \cite{xu2021duality}, in the case \( k = 2 \), \( r = 2p-1 \), we have (after adapting it to our current sign and summation order) that
	\begin{equation}\label{eqn:weighted}
		-T(\overline{2}, \{1\}_{2p-2}) = \sum_{j=1}^{2p-1} (-1)^{j-1} T(\overline{1}, \{1\}_{2p-2-j}) W(j+1, j) \,.
	\end{equation}
	Then from Equation (28) \cite{xu2021duality} (adapted to the current conventions), we have
	\begin{equation}\label{eqn:Tm111:eval}
		T(\overline{1}, \{1\}_{r}) = -\frac{(-1)^{\lfloor r/2 \rfloor}}{(r+1)!} \Big( \frac{\pi}{2} \Big)^{r+1} \,.
	\end{equation}
	Introduce the generating series
	\[
		G_{\bar1\{1\}}^T(z) \coloneqq \sum_{r=0}^\infty T(\overline{1}, \{1\}_r) z^r \,, \quad 
		\mathcal{W}(z) \coloneqq \sum_{j=1}^\infty W(j+1, j) z^j \,.
	\]
	From \autoref{eqn:Tm111:eval}, we readily have
	\begin{equation}\label{eqn:gs:Tm111}
		G_{\bar1\{1\}}^T(z) = \frac{1}{z} \Big( {-}1 + \cos\Big( \frac{\pi z}{2} \Big) - \sin\Big( \frac{\pi z}{2} \Big) \Big) \,.
	\end{equation}
	Taking \( \sum_{p=1}^\infty \bullet \) of both sides of \autoref{eqn:weighted},  and using \autoref{eqn:gs:Tm111}, and the evaluation in \autoref{lem:gs:Tm2111} gives the following lemma.
	\begin{Lem}\label{lem:W:full}
		The following generating series identity holds
		\begin{align*}
			 \frac{1}{2} \Big( \cos\Big( \frac{\pi z}{2} \Big) \big( \mathcal{W}(z) - \mathcal{W}(-z) \big) + \sin\Big( \frac{\pi z}{2} \Big) \big( \mathcal{W}(z) + \mathcal{W}(-z) \big) \Big)  
			 =  - 2 \cos\Big( \frac{\pi z}{2} \Big) D(z) \,.
		\end{align*}
	\end{Lem}
	\medskip
	
	Now we will evaluate \( W(2p+1, 2p) \) directly.  From Theorem 1.1 \cite{xu2021duality} in the case \( m = 1 \) (c.f. Equation (46) \cite{xu2021duality}), we have
	\begin{equation}\label{eqn:weighted:oddwt}
		W(2p+1, 2p) = \sum_{j=1}^p W(2j+1, 2) \cdot \frac{(-1)^{p-j} \pi^{2p-2j}}{(2p-2j+1)!} \,.
	\end{equation}
	(Again this expression has been adapted to the current sign conventions).  Here we have already substituted the factor \( Z(p,j) = \frac{(-1)^{p-j} \pi^{2p-2j}}{(2p-2j+1)!} \), c.f. Equations (9) and (10) in \cite{xu2021duality}.)
	On the other hand, using Theorem 3.4 \cite{xu2021duality} we have
	\begin{equation}\label{eqn:weighted:2}
		W(2k+1, 2) = 
		\frac{1}{2^{2k-1}} \sum_{j=1}^k (1 - 2^{2j+1}) (1 - 2^{1-(2k-2j)}) \zeta(2j+1) \zeta(2k-2j) \,.
	\end{equation}
	Taking \( \sum_{p=0}^\infty \bullet \) of both sides of \autoref{eqn:weighted:oddwt}, and using \autoref{eqn:weighted:2} (also in generating series form) straightforwardly gives the following lemma.
	\begin{Lem}\label{lem:W:odd}
		The following generating series identity holds
		\[
		\frac{1}{2} \big(\mathcal{W}(z) + \mathcal{W}(-z) \big) = \sum_{p=0}^\infty W(2p+1, 2p)  z^{2p} = \Big( A\Big( \frac{z}{2} \Big) - 2 A(z) \Big) \cos\Big( \frac{\pi z}{2} \Big) \,.
		\]
	\end{Lem}
	From \autoref{lem:W:full} and \autoref{lem:W:odd}, we obtain the following system of equations for \( \mathcal{W}(z) \), and \( \mathcal{W}(-z) \):
	\[
		\left\{ \, 
		\begin{aligned}[c]
		\begin{aligned}[b] \frac{1}{2} \Big( \!  \cos\Big( \frac{\pi z}{2} \Big) & \big( \mathcal{W}(z) - \mathcal{W}(-z) \big)  \\[-0.8ex]
		& {} + \sin\Big( \frac{\pi z}{2} \Big) \big( \mathcal{W}(z) + \mathcal{W}(-z) \big) \!  \Big) \end{aligned} &\,\, = \,\, - 2 \cos\Big( \frac{\pi z}{2} \Big) D(z)  \\[0.5ex]
		  \frac{1}{2} \big(\mathcal{W}(z) + \mathcal{W}(-z) \big) & \,\, = \,\, \Big( A\Big( \frac{z}{2} \Big) - 2 A(z) \Big) \cos\Big( \frac{\pi z}{2} \Big) \,.
		\end{aligned} \right.
	\]
	Solve this system of equations simultaneously for \( \mathcal{W}(z) \), and \( \mathcal{W}(-z) \), and we obtain the following proposition, through which we can evaluate \( W(j+1,j) \) in general.
	\begin{Prop}\label{prop:gs:W}
		The following generating series identity holds
		\[
			\mathcal{W}(z) = -2 D(z) + \Big( A\Big(\frac{z}{2}\Big) - 2A(z) \Big) \Big( \cos\Big( \frac{\pi z}{2} \Big) - \sin\Big( \frac{\pi z}{2} \Big) \Big)  \,.
		\]
	\end{Prop}

	Then from Equation (59) \cite{xu2021duality} in the case \( r = 2p, k = 2 \), (adapting things to our convention) we have that
	\begin{equation}\label{eqn:W:odd2}
		W(2p+1, 2p) = \sum_{j=1}^{2p} (-1)^{j-1} T(\overline{1}, \{1\}_{2p-j-1}) T(\overline{2}, \{1\}_{j-1}) \,.
	\end{equation}
	Taking \( \sum_{p=0}^\infty \bullet \) of both sides of \autoref{eqn:W:odd2} gives the following lemma.
	\begin{Lem}\label{lem:gs:W2}
		The following generating series identity holds
		\[
			\mathcal{W}(z) + \mathcal{W}(-z) = - G_{\bar2\{1\}}^T(z) \cdot \Big(\!  \cos\Big(\frac{\pi z}{2}\Big) - \sin\Big(\frac{\pi z}{2}\Big) \! \Big) - G_{\bar2\{1\}}^T(-z) \cdot \Big( \! \cos\Big(\frac{\pi z}{2}\Big) + \sin\Big(\frac{\pi z}{2}\Big) \! \Big)
		\]
	\end{Lem}
	
	Finally, from \autoref{lem:W:full} and \autoref{lem:gs:W2} we have the following system of equations for \( G_{\bar2\{1\}}(z) \), and \( G_{\bar2\{1\}}(-z) \):
	\[
		\left\{ \,\,
		\begin{aligned}[c]
		& 
		 \frac{1}{2} \Big( G_{\bar2\{1\}}^T(z) - G_{\bar2\{1\}}^T(-z) \Big) \,\, = \,\,
		\begin{aligned}[t] \frac{1}{2} \Big( \! \cos\Big( \frac{\pi z}{2} \Big) & \big( \mathcal{W}(z) - \mathcal{W}(-z) \big) \\[-0.8ex]
			&  {} + \sin\Big( \frac{\pi z}{2} \Big) \big( \mathcal{W}(z) + \mathcal{W}(-z) \big) \!  \Big)  \end{aligned} \\
		& \begin{aligned}[b] 
			- G_{\bar2\{1\}}^T(z) \cdot & \Big( \!  \cos\Big(\frac{\pi z}{2}\Big) - \sin\Big(\frac{\pi z}{2}\Big) \! \Big) \\[-0.5ex]
			& - G_{\bar2\{1\}}^T(-z) \cdot \Big( \!  \cos\Big(\frac{\pi z}{2}\Big) + \sin\Big(\frac{\pi z}{2}\Big)\!  \Big) \end{aligned} \,\, =  \,\, \mathcal{W}(z) + \mathcal{W}(-z) \,. \\[1ex] 
		\end{aligned}
		\right.
	\]
	Solve this system of equations simultaneously for \( \smash{G_{\bar2\{1\}}^T(z)} \), and \( \smash{G_{\bar2\{1\}}^T(-z)} \), and eliminate \( \mathcal{W}(z) \) via \autoref{prop:gs:W}.  We obtain the following proposition, evaluating \( T(\overline{2}, \{1\}_{r-1}) \) in general.
	\begin{Prop}\label{prop:gs:Tm2111}
		The following generating series identity holds
		\[
			G_{\bar2\{1\}}^T(z) = -A\Big( \frac{z}{2} \Big) + 2A(z) - 2D(z) \Big( \!  \cos\Big( \frac{\pi z}{2} \Big) + \sin\Big( \frac{\pi z}{2} \Big) \! \Big) \,.
		\]
	\end{Prop}

	Recall also, the duality relation
			\begin{equation*}
			 T(\overline{2}, \{1\}_m) = - (-1)^m T(\overline{1}, \{1\}_m, \overline{1}) \,,
			 \end{equation*}
	from \cite[Theorem 4.3]{xu2021duality}, in the case \( p = 1, r = m+1 \).  (The expression in \cite{xu2021duality} has again been adapted to the sign and order conventions of the current paper.)  We therefore also obtain an evaluation for \( T(\overline{1}, \{1\}_m, \overline{1}) \) from this proposition. 

	\paragraph{First expression for \( Q(z) \)} Substituting the result of \autoref{prop:gs:Tm2111} into \autoref{lem:Q:expr1}, we obtain the following expression for the generating series \( Q(z) \).  (Some simplification via \( A(-z) = A(z) \), and \( D(-z) = -D(z) \) is necessary.)
	\begin{Prop}[First expression for \( Q(z) \)]\label{prop:Q:eval1}
		The following generating series evaluation holds
	\begin{equation}\label{eqn:Q:eval1}
		Q(z) = \frac{z}{2} \exp\Big( \frac{\pi z}{4} \Big) \Big(  2 D'\Big(\frac{\ii z}{2} \Big) - \ii\pi D\Big(\frac{\ii z}{2}\Big) \! \Big) - \frac{\ii z}{4} \Big( A'\Big(\frac{\ii z}{4} \Big) - 4 A'\Big( \frac{\ii z}{2} \Big) \! \Big) \,.
	\end{equation}
	\end{Prop} \medskip

	Continuing with Xu, Yan and Zhao's argument: In order to derive the second expression for \( Q(z) \), we first recall the evaluation of certain cotangent integrals, and use this to derive generating series for arctangent integrals.  Then Theorem 4.5 \cite{xu2022alternating} allows us to express \( Q(z) \) via \( T(2, \{1\}_{2p}) \), \( S(2, \{1\}_{2p+1}) \), these arctangent integral, and the new values of interest \( T(2, \{1\}_{2p-1}, \overline{1}) \) and \( S(2, \{1\}_{2p-2}, \overline{1}) \). 

	\paragraph{Cotangent integrals} We continue by recalling the following identities from Lemma 4.1 of \cite{xu2022alternating} (given originally in \cite[Equations (2.2) and (2.3)]{orr2017generalized}).  Write \( \delta_\bullet \) to be \( 1 \) if the condition in \( \bullet \) is true, and \( 0 \) otherwise (also write \( \lfloor \bullet \rfloor \) rather than \( [\bullet] \) for clarity).  For \( p \geq 1 \in \mathbb{Z} \), we have
	\begin{align}
	\label{eqn:pi2:cot}
	\int_0^{\pi/2} x^p \cot(x)dx&=
	\begin{aligned}[t]
	& \left(\frac{\pi}{2}\right)^p \bigg\{ \log(2) + \sum_{k=1}^{\lfloor p/2 \rfloor} {}  \frac{p!(-1)^k(4^k-1)}{(p-2k)!(2\pi)^{2k}}\zeta(2k+1) \bigg\} \\[-0.5ex]
	&{} +\delta_{\text{$p$ even}} \, \frac{p!(-1)^{p/2}}{2^p}\zeta(p+1) \,,
	\end{aligned}\\[1ex]
	\label{eqn:pi4:cot}
	\int_0^{\pi/4} x^p \cot(x)dx&=
	\begin{aligned}[t]
	&\frac{1}{2}\left(\frac{\pi}{4}\right)^p \bigg\{ 
	\begin{aligned}[t]  
	{} \log(2) + {} & \sum_{k=1}^{\lfloor p/2 \rfloor}  \frac{p!(-1)^k(4^k-1)}{(p-2k)!(2\pi)^{2k}}\zeta(2k+1) \\[-0.5ex]
	&  -\sum_{k=1}^{\lfloor (p+1)/2 \rfloor} \frac{p!(-4)^k\beta(2k)}{(p+1-2k)!\pi^{2k-1}} \bigg\} \\[-0.5ex]
	\end{aligned} \\
	&\quad+\delta_{\text{$p$ even}} \, \frac{p!(-1)^{p/2}}{2^p}\zeta(p+1) \,,
	\end{aligned}
	\end{align}
	Introduce the following generating series of these cotangent integrals
	\begin{align*}
	C_{\half}(z) \coloneqq \sum_{p=1}^\infty \frac{z^p}{p!} \int_{0}^{\pi/2} x^p \cot(x) \,\mathrm{d}x \,, \quad 
	C_{\quart}(z) \coloneqq \sum_{p=1}^\infty \frac{z^p}{p!} \int_{0}^{\pi/4} x^p \cot(x) \,\mathrm{d}x \,.
	\end{align*}
	Taking \( \sum_{p=1}^\infty \frac{z^p}{p!} \bullet \) of both sides of Equations \eqref{eqn:pi2:cot} and \eqref{eqn:pi4:cot}  leads straightforwardly to the following lemma.
	
	\begin{Lem}
		The following generating series evaluations hold
		\begin{align*}
		C_{\half}(z) &= A\Big(\frac{\ii z}{2}\Big) -\log (2) + \exp\Big( \frac{\pi  z}{2} \Big) \Big({-}A\Big(\frac{\ii z}{4}\Big)+A\Big( \frac{\ii z}{2}\Big)+\log (2)\Big) \,, \\ 
		C_{\quart}(z) &= A\Big(\frac{\ii z}{2}\Big) - \frac{1}{2} \log (2) + \frac{1}{2}  \exp\Big(\frac{\pi  z}{4}\Big) \Big({-}A\Big(\frac{\ii z}{8}\Big) + A\Big(\frac{\ii z}{4} \Big) - 2 \ii D\Big(\frac{\ii z}{2}\Big)+\log (2) \Big)  \,.
		\end{align*}
	\end{Lem}
	
	\paragraph{Arctangent integrals} Now introduce the generating series of the arctangent integral
	\[
	R(z) \coloneqq \sum_{p=1}^\infty \frac{z^p}{p!} \int_0^1 \arctan(x)^p \, \mathrm{d}x \,.
	\]
	From the proof of Proposition 4.2 (Equation (4.3) specifically) in \cite{xu2022alternating}, we have --  after substituting \( x = \tan(t) \) and integrating by parts -- that, for \( p \geq 1 \in \mathbb{Z} \), 
	\begin{align*}
	& \int_0^1 \arctan^p(x) \, \mathrm{d}x = 
	\Big(\frac{\pi}{4}\Big)^p - \frac{p \, \pi^{p-1}}{2^p} - p \sum_{k=1}^{p-1} (-1)^k \binom{p-1}{k} \Big( \frac{\pi}{2} \Big)^{p-1-k} \int_{\pi/4}^{\pi/2} u^k \cot(u) \, \mathrm{d}u \,.
	\end{align*}
	Taking \( \sum_{p=1}^\infty \frac{z^p}{p!} \bullet \) of both sides of this leads immediately to the following lemma.  (Here we have also used that \( A(-z) = A(z) \), and \( D(-z) = -D(z) \).)
	\begin{Lem}\label{lem:R}
		The following generating series identities hold
		\begin{align*}
		R(z) &= -1 + \exp\Big( \frac{\pi z}{4} \Big) - z \exp\Big( \frac{\pi z}{2} \Big) \cdot \Big( \frac{1}{2} \log(2) - C_{\half}(-z) - C_{\quart}(-z) \Big) \\[0.5ex]
		& = \begin{aligned}[t]
		-1 + \exp\Big( \frac{\pi z}{4} \Big) & - z \Big({-}A\Big(\frac{\ii z}{4}\Big)+A\Big( \frac{\ii z}{2}\Big)+\log (2)\Big) \\[-0.5ex]
		& {} +  \frac{z}{2} \exp\Big( \frac{\pi z}{4} \Big) \Big({-}A\Big(\frac{\ii z}{8}\Big) + A\Big(\frac{\ii z}{4} \Big) + 2 \ii D\Big(\frac{\ii z}{2}\Big)+\log (2) \Big) \,.
		\end{aligned}
		\end{align*}
	\end{Lem}
	
	\paragraph{The ``arctangent over $x$'' integral}  In Theorem 4.5 of \cite{xu2022alternating}, the authors establish some (rather complicated looking) expressions for the ``arctangent over $x$'' integral.  The original expression uses the alternating \( t \) value, \( t(\overline{k}) \coloneqq -\beta(k) \); we have written this directly using \( t(\overline{1}) = \frac{\pi}{4} \), and \( t(\overline{2}) = -\beta(2) = -G \), with \( G \) the Catalan constant, likewise \( t(2) \coloneqq \sum_{n=1}^\infty \frac{1}{(2n-1)^2} = \frac{\pi^2}{8} \) is a standard result.  We also need to introduce \( r(p) \coloneqq \int_0^1 \arctan^p(x) \, \mathrm{d}x \).  (This is denoted by \( A(p) \) in \cite{xu2022alternating}, but this clashes with the notation \( A(z) \) for the zeta generating series, taken from \cite{zagier2012evaluation}.)
	{
	\begin{align}
	\label{eqn:arctanoverx:even}
	\begin{aligned}
 	\mathllap{\int_0^1 \frac{\arctan^{2m}(x)}{x} \, \mathrm{d}x
	= {}}
	 &\Big(\frac{\pi}{4}\Big)^{2m-1}\bigg({-}\beta(2)+\frac{\pi^2}{8}\bigg)+\frac{(-1)^{m}(2m-1)!}{2^{2m}} \, S(2,\{1\}_{2m-2},\overline{1})\\
	\hspace{8em} &{}+(2m-1)!\sum_{u=0}^{m-1}\frac{(-1)^u \, r(2m-2u-1)}{(2m-2u-1)! \,  2^{2u+1}} \, T(2,\{1\}_{2u})\\
	&{}-(2m-1)!\sum_{v=1}^{m-1}\frac{(-1)^{v+1} \, (\pi/4)^{2m-2v-1}}{2^{2v+1} \,  (2m-2v-1)!} \Big( T(2,\{1\}_{2v})+T(2,\{1\}_{2v-1},\overline{1}) \Big)\\
	&{}-(2m-1)!\sum_{v=1}^{m-1}\frac{(-1)^{v} \, (\pi/4)^{2m-2v}}{2^{2v} \,  (2m-2v)!} \Big( S(2,\{1\}_{2v-1})-S(2,\{1\}_{2v-2},\overline{1})\Big)
	\end{aligned}
	\\[1ex]
	\label{eqn:arctanoverx:odd}
	\begin{aligned}
	\mathllap{\int_0^1 \frac{\arctan^{2m+1}(x)}{x} \, \mathrm{d}x 
	= {}}
	&{}-\Big(\frac{\pi}{4}\Big)^{2m} \bigg( {-}\beta(2) + \frac{\pi^2}{8}\bigg) - \frac{(-1)^{m}(2m)!}{2^{2m+1}} \,  T(2,\{1\}_{2m-1},\overline{1})\\
	&{}+(2m)!\sum_{u=0}^{m-1}\frac{(-1)^u \, r(2m-2u)}{(2m-2u)! \,  2^{2u+1}} \,  T(2,\{1\}_{2u})\\
	&{}-(2m)!\sum_{v=1}^{m-1}\frac{(-1)^{v} \, (\pi/4)^{2m-2v}}{2^{2v+1} \, (2m-2v)!} \Big(T(2,\{1\}_{2v})+T(2,\{1\}_{2v-1},\overline{1}) \Big)\\
	&{}-(2m)!\sum_{v=0}^{m-1}\frac{(-1)^{v} \, (\pi/4)^{2m-2v-1}}{2^{2v+2} \,  (2m-2v-1)!} \Big(S(2,\{1\}_{2v+1})-S(2,\{1\}_{2v},\overline{1}) \Big),
	\end{aligned}
	\end{align}
	}
	
	If we introduce the final few generating series (\( E \) and \( F \) are the important ones from the introduction)
	\begin{alignat*}{4}
		G_{2\{1\}}^S(z) & \coloneqq \sum_{m=0}^\infty S(2, \{1\}_{2m+1}) z^{2m+2} \,, & \quad\quad 
		E(z) &\coloneqq  \sum_{m=0}^\infty S(2, \{1\}_{2m}, \overline{1}) z^{2m+2} \,,  \\
		G_{2\{1\}}^T(z) & \coloneqq \sum_{m=0}^\infty T(2, \{1\}_{2m}) z^{2m+1} \,, &
		F(z) & \coloneqq \sum_{m=0}^\infty T(2, \{1\}_{2m+1}, \overline{1}) z^{2m+3} \,, 
	\end{alignat*}
	then by taking 
	\[ 
	\sum_{m=1}^\infty \frac{(-2 \ii z)^{2m}}{(2m-1)!} \cdot \text{(Eqn.~\ref{eqn:arctanoverx:even})} \quad + \quad \sum_{m=1}^\infty \frac{(-2 \ii z)^{2m+1}}{(2m)!} \cdot \text{(Eqn.~\ref{eqn:arctanoverx:odd})} \,,
	\]
	we obtain (after some work) the following alternative expression for \( Q(-2\ii z) \).  (Note the choice \( -2\ii z \), in order to makes the final deduction below slightly more direct.)
	\begin{Prop}[Second expression for \( Q(z) \)]\label{prop:Q:eval2}
		The following generating series identity holds
		\begin{equation}
		\label{eqn:Q:eval2}
			 \begin{aligned}[c]
			 Q(-2\ii z) = {}
				& -2\ii \beta(2) z + \Big( 1 - \exp\Big( {-}\frac{\ii\pi z}{2} \Big) \! \Big) \big( 2 \ii \beta(2) z + G_{2\{1\}}^S(z) - \ii G_{2\{1\}}^T(z) \big) \\
				& - \ii R(-2\ii z) G_{2\{1\}}^T(z) + \exp\Big( {-}\frac{\ii \pi z}{2} \Big) \big( E(z) + \ii  F(z) \big) \,.
			\end{aligned}
		\end{equation}
	\end{Prop}
	
	\paragraph{Two remaining evaluations} Before we can equate the two expressions for \( Q(z) \), or rather \( Q(-2\ii z) \), and extract a useful result for \( E(z) + \ii F(z) \), we still need to evaluate \( S(2, \{1\}_{2p-1}) \) and \( T(2, \{1\}_{2p}) \).  Fortunately, this is straightforward.
	
	The multiple \( T \) values satisfy the same duality as multiple zeta values (see \cite[\S3.1]{kaneko2020level2}).  In particular 
	\begin{equation}\label{eqn:T2111}
	 T(2, \{1\}_r) = T(r+2) = 2(1 - 2^{-2r-2}) \zeta(2r + 2) \,.
	 \end{equation}
	 Whereas from \cite[Equation (3.17)]{xu2022variants}, after rewriting the alternating zeta values \( \overline{\zeta}(k) \coloneqq -\zeta(\overline{k}) = (1 - 2^{1-k}) \zeta(k) \) directly in terms of normal zeta values, except for \( \overline{\zeta}(1) \coloneqq -\zeta(\overline{1}) = \log(2) \)  (and adapting to our conventions), we have
	 \begin{equation}\label{eqn:S2111}
	 	S(2, \{1\}_{2p-1}) = \begin{aligned}[t] 
	 	& 2p T(2p+1) - 2 \log(2) T(2p) \\[-0.5ex]
	 	& - 2 \sum_{j=0}^{p-2} (1 - 2^{1-(2p-1-2j)}) \zeta(2p-1-2j) T(2j+2) 
	 	\end{aligned}
	 \end{equation}
	 From \autoref{eqn:T2111} and \autoref{eqn:S2111} we readily obtain the final lemma.
	 \begin{Lem}\label{lem:ST2111}
	 	The following generating series evaluations hold
	 	\begin{align*}
	 		G_{2\{1\}}^T(z) &= \frac{\pi}{2} \tan\Big( \frac{\pi z}{2} \Big) \,, \\
	 		G_{2\{1\}}^S(z) &= \frac{\pi z}{2} \tan\Big( \frac{\pi z}{2} \Big) \Big( \! A\Big( \frac{z}{2} \Big) - A(z) - \log(2) \Big) - \frac{z}{2} A'\Big( \frac{z}{2} \Big) + 2 z A'(z) \,,
	 	\end{align*}
	 \end{Lem}
	 
	\paragraph{Conclusion} Substitute the results of \autoref{lem:ST2111}, and \autoref{lem:R} into \autoref{prop:Q:eval2}, and rewrite \( Q(-2\ii z) \) using the earlier expression from \autoref{prop:Q:eval1}.  This gives an identity involving only \( E(z) + \ii F(z) \) and known generating series in terms of \( A, D \) and trigonometric/exponential functions.  Solving this identity for \( E(z) + \ii F(z) \), and simplifying the result gives \autoref{thm:qn1:eval}. \hfill \qedsymbol
	
	\section{Recap of Goncharov's iterated integrals} \label{sec:goncharov}
	
	In order to tackle the claims of Question 2 \cite[p.~21]{xu2022alternating}, and establish the results in \autoref{thm:qn2:mtv} and \autoref{thm:qn2:msv}, we need to recall some definitions, properties and results of the Chen iterated integrals \cite{chen77}, and multiple polylogarithms.  We refer mainly to the theory developed by Goncahrov \cite[\S2]{goncharov01}, although we convert the results therein to the current \( > \)-convention for MPL's and MZV's used in the present article.
	
	\paragraph{Definitions} For a family of differential one-forms \( \omega_1, \ldots, \omega_n \) on a manifold \( M \), and a path \( \gamma \colon [0,1] \to M \) (piecewise-)smooth, the iterated integral is inductively defined \cite[c.f. Eqn.~(5)]{goncharov01} as
	\[
		\int_\gamma \omega_1 \circ \cdots \circ \omega_n \coloneqq \int_0^1 \gamma^\star \omega_1 \Big( \int_{\gamma\rvert_{[0,t]}} \omega_2 \circ \cdots \circ \omega_n \Big) \,.
	\]
	Then we define the iterated integral function (or hyperlogarithm) as
	\begin{equation}\label{eqn:II:def}
		I(a; x_1,\ldots,x_n; b) \coloneqq \int_b^a \frac{\mathrm{d}t}{t-x_1} \circ \cdots \circ \frac{\mathrm{d}t}{t-x_n} 
		= \int_{a > t_1  > \cdots > t_n > b} \frac{\mathrm{d}t_1}{t_1 - x_1} \cdots \frac{\mathrm{d}t_n}{t_n - x_n} \,,
	\end{equation}
	taken along the straight-line path from \( b \) to \( a \).  (This, perhaps unusual, convention gives the integration region \( a > t_1 > \cdots > t_n > b \), which aligns more directly with the notation in \cite{xu2022alternating}.)  This is convergent when \( x_1 \neq a \) and \( x_n \neq b \), as the logarithmic singularity still leads to convergent integrals.
	
	The convergent integrals satisfy an affine invariance; when \( a, b, x_i \) are all simultaneously transformed under \( f(x) = \alpha x + \beta \), we have
	\begin{equation}\label{eqn:affine}
		I(a; x_1,\ldots,x_n; b) = I(f(a); f(x_1), \ldots, f(x_n); f(b)) \,.
	\end{equation}
	This is directly obtained by the change of variables \( x_i \mapsto f(x_i) \) in \autoref{eqn:II:def}.
	
	\paragraph{Derivative and differentials} Then Goncharov \cite[Theorem 2.1]{goncharov01} establishes the following differential equation satisfies by the iterated integral \( I \).  Namely
	\begin{equation}\label{eqn:II:totaldiff}
	\begin{aligned}
		& \mathrm{d}I(x_0; x_1,\ldots,x_n; x_{n+1})  \\
		& = \sum\nolimits_{i=1}^n I(x_0; x_1,\ldots,\widehat{x_i},\ldots,x_n; x_{n+1}) \cdot  \mathrm{d} \bigl( \log( x_i - x_{i-1}) - \log( x_i - x_{i+1} )  \bigr)  \,,
	\end{aligned}
	\end{equation}
	where \( \widehat{x_i} \) denotes that \( x_i \) is dropped from the argument string.
	This follows directly by passing the derivative through the integral, and computing the result via partial fractions identity.  (See the proof of Theorem 2.1 \cite{goncharov01} for details.)
	
	The proof of this differential equation leads to the following very useful result: for \( a_i \) constants with respect to \( y \),
	\begin{equation}\label{eqn:II:diff}
		\frac{\mathrm{d}}{\mathrm{d}y} I(y; a_1,\ldots,a_m; a_{m+1}) = I(y; a_2,\ldots,a_m; a_{m+1}) \cdot \frac{1}{y - a_1} \,,
	\end{equation}
	so that a primitive with respect to \( y \) of the right-hand side is given by inserting \( a_1 \) at the start of that iterated integral. 
	
	\paragraph{Shuffle product} We also need to recall that the iterated integrals can be multiplied using the shuffle-product to obtain an algebra structure \cite[c.f. Eqn.~(1.5.1)]{chen71} (wherein it is attributed as an observation due to Ree \cite{ree58}).  Namely 
	\[
		\int_{\gamma} \omega_1 \circ \cdots \circ \omega_r \cdot \int_{\gamma} \omega_{r+1} \circ \cdots \circ \omega_{r+s} = \sum_{\sigma \in \Sigma(r,s)} \int_{\gamma} \omega_{\sigma(1)} \circ \cdots \circ \omega_{\sigma(r+s)} \,,
	\]
	where
	\[
		\Sigma(r,s) = \bigl\{ \sigma \in S_{r+s} \mid \sigma^{-1}(1) < \cdots < \sigma^{-1}(r) \,, \sigma^{-1}(r+1) < \cdots < \sigma^{-1}(r+s) \bigr\} \,,
	\]
	is the set of so-called \( (r,s) \)-shuffles.  In particular, the same result holds for \( I(b; x_1,\ldots,x_r; a) \cdot I(b; x_{r+1}, \ldots, x_{r+s}; a) \), where all \( (r,s) \)-shuffles of the \( x_i \) are taken.
	
	This is a manifestation of writing all compatible ways of interleaving the sets of inequalities \( a > t_1 > \cdots > t_r > b \) from the first integral and \( a > t'_1 > \cdots > t'_s > b \) from the second integral, when writing the product as an \( (r+s) \)-dimensional integral via Fubini.  Therein the case \( t_i = t'_j \) defines a set of measure 0, and so can be neglected.

	\paragraph{Regularisation of integrals} The `canonical' asymptotic regularisation \cite[c.f. \S2.9]{goncharov01}, given by writing
	\[
	I(x-\eps; 0; \eps) = \int_{\eps}^{x-\eps} \frac{\mathrm{d}t}{t} = \log(x + \eps) - \log(\eps) \,,
	\]
	as a polynomial in \( \log(\eps) \) with convergent coefficients (as $\eps\to0$), allowed one to formally the notion of iterated integrals to the case \( x_1 = a \) or \( x_n = b \).  One obtains the shuffle-product induced regularisation which sets \( I(a;b;c) = \log(b-a) - \log(b-c) \), where \( \log(0) \coloneqq 0 \).  See \S2.9 \cite{goncharov01}, and the shuffle product introduced below.  Compare also with the remark before Eqn.~(22) \cite{goncharov01}.  This allows one to extend the differential formulae, \autoref{eqn:II:totaldiff} and \eqref{eqn:II:diff} above, to the case where \( x_i = x_{i\pm 1} \), by regularising \( \log(x_i - x_{i\pm1}) = 0 \).

	\paragraph{Relation to multiple polylogarithms}  Recall the multiple polylogarithms are defined by the following conical sum
	\[
		\Li_{s_1,\ldots,s_r}(x_1,\ldots,x_r) \coloneqq \sum_{m_1 > \cdots > m_r > 0} \frac{x_1^{m_1} \cdots x_r^{m_r}}{m_1^{s_1} \cdots m_r^{s_r}} \,, \quad |x_i| < 1 \,.
	\]
	In this notation the alternating multiple zeta values of \autoref{eqn:zeta} are just 
	\[
	\zeta^{\sigma_1,\ldots,\sigma_r}(s_1,\ldots,s_r) = \Li_{s_1,\ldots,s_r}(\sigma_1,\ldots,\sigma_r)  \,,
	\] and the alternating multiple mixed values \( M_{\sigma_1,\ldots,\sigma_r}^{\eps_1,\ldots,\eps_r}(s_1,\ldots,s_r) \) can be expressed via suitable \( \mathbb{Q}[\ii] \)-linear combinations of \( \Li_{s_1,\ldots,s_r}(x_1,\ldots,x_r) \), \( x_1 \in \{ \pm 1, \pm \ii \} \), after expanding out the numerator of \autoref{eqn:MM:def}.
	
	By expanding out the geometric series, and term-by-term integration (c.f. Theorem 2.2 \cite{goncharov01}), one can straightforwardly show (for \( s_i \in \mathbb{Z}_{>0} \)) that
	\begin{equation}\label{eqn:mpl}
		\Li_{s_1,\ldots,s_r}(x_1,\ldots,x_r) = (-1)^r I(1; \{0\}_{s_1-1}, \tfrac{1}{x_1}, \{0\}^{s_2-1}, \tfrac{1}{x_1 x_2} , \ldots, \{0\}^{s_r-1}, \tfrac{1}{x_1 x_2 \cdots x_r}; 0) \,,
	\end{equation}
	in terms of the iterated integrals. \medskip
	
	\noindent We are now in a position to apply these objects and properties to the evaluation of \( T(\overline{2}, 1, \ldots, 1, \overline{1}) \) and \( S(\overline{2}, 1,\ldots,1,\stackon[.1pt]{$1$}{\brabar}) \) which arise in \autoref{thm:qn2:mtv} and \autoref{thm:qn2:msv}.
	
	\section{Proof of \autoref{thm:qn2:mtv}} \label{sec:thm2a}
	
	The iterated integral representation of MTV's \cite[\S1.2]{xu2022alternating} allows us to write
	\[
		T(\overline{2}, \{1\}_{m-1}, \overline{1}) = (-1)^{\lfloor (m+1)/2 \rfloor} \int_0^1 \frac{\mathrm{d}t}{t} \circ \Big\{ \frac{-2 \,  \mathrm{d}t}{t^2 + 1} \Big\}_m \, \circ  \frac{-2 \, \mathrm{d}t}{t^2 - 1}  \,.
	\]
	Inside an iterated integral, the notation \( \{ f(t)\mathrm{d}t \}_m \) should be interpreted as
	\[
		\underbrace{f(t) \mathrm{d}t \circ \cdots \circ f(t) \mathrm{d}t}_{\text{$m$ repetitions}} \,.
	\]
	We prefer \( \{ f(t) \mathrm{d}t \}_m \), instead of just \( (f(t)\mathrm{d}t)^m \), to avoid possible confusion with \( f(t)^m \mathrm{d}t \); the notation immediately clarifies this is an iterated integral, not a one-dimensional integral.   Then since
	\begin{equation}\label{eqn:arctanint}
		\int_r^s \Big\{ \frac{-2 \,\mathrm{d}t}{t^2 + 1} \Big\}_m = \frac{1}{m!} \Big( \int_r^s \frac{-2 \,  \mathrm{d}t}{t^2 + 1} \Big)^m = \frac{(-2)^m}{m!} (\arctan(s) - \arctan(r))^m \,,
	\end{equation}
	we therefore have
	\[
		T(\overline{2}, \{1\}_{m-1}, \overline{1}) = \frac{1}{m!} (-1)^{\lfloor (m+1)/2 \rfloor} \int_{1 > s > r > 0}  (-2)^m (\arctan(s) - \arctan(r))^m \cdot \frac{\mathrm{d}s}{s} \cdot \frac{-2 \, \mathrm{d}r}{r^2 - 1} \,,
	\]
	as a general iterated integral of length 2.  Make the substitution
	\begin{equation}\label{eqn:int:subst}
		(r,s) = \big( \tan(-\tfrac{1}{2\ii} \log(x)), \tan(-\tfrac{1}{2\ii} \log(y)) \big) \,,
	\end{equation}
	and we obtain
	\[
		T(\overline{2}, \{1\}_{m-1}, \overline{1}) = \frac{\ii}{m!} (-1)^{\lfloor (m+1)/2 \rfloor} \cdot \ii^m \int_{\gamma \,,\, -\ii > y > x > 1} \frac{4 (\log(x) - \log(y))^m}{(x^2 + 1)(y^2 - 1)} \mathrm{d}x \, \mathrm{d}y  \,,
	\]
	clockwise along the circular arc \( \gamma \) from \( 1 \) to \( -\ii \).  (This integral is homotopy invariant, with no singularities (in the respective variables) at the end points so we can deform the path to the straight line, say.) The goal is now to understand and evaluate
	\[
		I_{a,b} \coloneqq \int_{-\ii>y>x>1} \frac{4}{a! \, b!} \frac{\log(x)^a \log(y)^b}{(x^2+1)(y^2 - 1)} \, \mathrm{d}x \, \mathrm{d}y \,,
	\]
	which we can then piece together to via the binomial theorem to evaluate \( T(\overline{2}, \{1\}_{m-1}, \overline{1}) \).  Overall, after noticing \( \ii \cdot (-1)^{\lfloor (m+1)/2 \rfloor} \cdot \ii^m  = \ii^{\delta_{\text{$m$\,even}}} \), we obtain
	\begin{equation}\label{eqn:tm2111m1:asint}
		T(\overline{2}, \{1\}_{m-1}, \overline{1}) 
		= \ii^{\delta_{\text{$m$\,even}}} \sum_{\substack{a+b=m \\ a,b \geq 0}} (-1)^b I_{a,b} \,.
	\end{equation}
	(Note: \( I_{a,b} \) includes the factorials from the denominator of \( \binom{m}{a} = \frac{m!}{a! \, b!} \), hence the simple form.) 
	
	\paragraph{Primitive with respect to \( y \)} As a first step, we can directly write down a primitive for the integration with respect to \( y \) appearing in \( I_{a,b} \).  Under the shuffle-regularisation prescription (see \autoref{sec:goncharov} above)
	\[
		\frac{1}{b!} \log(y)^b = I(y; \{0\}_b; 0) \,,
	\]
	and since
	\[
		\frac{2}{y^2 - 1} =  \frac{1}{y - 1} - \frac{1}{y - (-1)}  \,,
	\]
	we immediately have (using the differential behaviour of the iterated integrals explained in \autoref{eqn:II:diff}) that
	\[
		\int \frac{2}{b!} \frac{ \log(y)^b}{y^2 - 1} \mathrm{d}y =  I(y; 1, \{0\}_b; 0) -  I(y; -1, \{0\}_b; 0) \,.
	\]
	
	So on taking \( [\bullet]_{y=x}^{y=-\ii} \) of the above primitive, we obtain
	\[
		\int_{-\ii > y > x} \frac{2}{b!} \frac{ \log(y)^b}{y^2 - 1} \mathrm{d}y \, = \,
		\begin{aligned}[t]
				& I(-\ii; 1, \{0\}_b; 0) -  I(-\ii; -1, \{0\}_b; 0) \\
				& \hspace{3em} {} - I(x; 1, \{0\}_b; 0) +  I(x; -1, \{0\}_b; 0)
			 \,.
		\end{aligned}
	\]
	
	\paragraph{Primitive with respect to \( x \)}  To continue with the \( x \) integration, we can first write the integrand as a  ($\mathbb{C}$-)linear combination of integrals whose upper bound is \( x \), and apply the same routine as before.  In particular, by the shuffle product of iterated integrals, we have
	\begin{align*}
		\frac{1}{a!} \log^a(x) \cdot I(x; \pm1, \{0\}_b; 0) 
		&= I(x; \{0\}_a; 0) \cdot I(x; \pm1, \{0\}_b; 0) \\
		&= \sum_{\substack{p+q=a \\ p,q \geq 0}} \binom{b+p}{p} I(x; \{0\}_q, \pm 1, \{0\}_{b+p}; 0) \,.
	\end{align*}
	(Note that since \( I(-\ii; \pm 1, \{0\}_b; 0) \) is just some complex number constant, we can leave the term \( \frac{1}{a!} \log^a(x) \cdot I(-\ii; \pm 1, \{0\}_b; 0) \) as it is in the integrand, and simply pull the constant out of the integral.)  Then since
	\[
		\frac{2}{x^2 + 1} = -\ii \bigg( \frac{1}{x - \ii} - \frac{1}{x - (-\ii)} \bigg) \,,
	\]
	we can straightforwardly obtain (again using \autoref{eqn:II:diff}) that
	\begin{align}
		\nonumber I_{a,b} & = \int_{-\ii > x > 1} \frac{2}{a!} \frac{\log(x)^a}{x^2 + 1} \cdot \begin{aligned}[t]
		\big( & I(-\ii; 1, \{0\}_b; 0) -  I(-\ii; -1, \{0\}_b; 0) \\[-0.5ex]
		& \hspace{3em} {} - I(x; 1, \{0\}_b; 0) +  I(x; -1, \{0\}_b; 0) \big) \, \mathrm{d}x 
		\end{aligned} \\
		\label{eqn:Iab:part1}& = 
			 -\ii \bigg[ 
			 \begin{aligned}[t] 
			& \,\, \big( I(x; \ii, \{0\}_a, 0) - I(x; -\ii; \{0\}_a; 0)\big) 
			\cdot \big( I(-\ii; 1, \{0\}_b; 0) - I(-\ii; -1, \{0\}_b; 0) \big) 
			 \end{aligned} 
			 \\[1.0ex]
			 &  \label{eqn:Iab:part2} \begin{aligned}[c]
			\smash{\phantom{ = -\ii \bigg[ } +\sum_{\substack{p+q = a \\ p,q \geq 0}} \binom{b+p}{p}} \big( 
			& -  I(x; \ii, \{0\}_q, 1, \{0\}_{b+p}; 0)
			+ I(x; \ii, \{0\}_q, -1, \{0\}_{b+p}; 0) \\
			& + I(x; -\ii, \{0\}_q, 1, \{0\}_{b+p}; 0)
			- I(x; -\ii, \{0\}_q, -1, \{0\}_{b+p}; 0)
			\big) \,\, \smash{\bigg]_{x=1 \,.}^{x=-\ii}}  \end{aligned}
	\end{align}
	It will be useful to write \( \widetilde{I_{a,b}} = \widetilde{I_{a,b}}(x) \) for the \( x \)-primitive before  any substitution, so we can analyse the behaviour at \( x = -\ii \) and \( x = 1 \) more carefully.  That is 
	\[ \widetilde{I_{a,b}} = -\ii \big[ \text{$(I+I)\cdot(I+I)$ terms in Eqn.~\eqref{eqn:Iab:part1}} 
	+ \text{summation in Eqn.~\eqref{eqn:Iab:part2}} \big] \,,
	\]
	so that
	\[
		I_{a,b} = \big[ \widetilde{I_{a,b}}(x) \big]_{x=1 \,.}^{x=-\ii}
	\]
		
	\paragraph{Simplification when summing \( (-1)^b \stackon[-4pt]{$I_{a,b}$}{\bratw} \)}  Firstly, let us note a large simplification in the binomial sum above, which occurs after computing
	\begin{equation}\label{eqn:Iab:sumreg}
		\sum_{\substack{a+b=m \\ a,b \geq 0}} (-1)^b \bigg( \sum_{\substack{p+q=a \\ p,q \geq 0}} \binom{b+p}{p} I(C; B, \{0\}_{q}, A, \{0\}_{b+p}; 0) \bigg) \,,
	\end{equation}
	in general.  We find that this is
	\begin{align*}
			& = \sum_{\substack{p+q+b=m \\ p,q,b \geq 0}} (-1)^b \binom{b+p}{p} I(C; B, \{0\}_{q}, A, \{0\}_{b+p}; 0) \\
			& = \sum_{\substack{N + q = m \\ N,q \geq 0}} \sum_{\substack{b+p = N \\ b,p \geq 0}} (-1)^b \binom{N}{p} I(C; B, \{0\}_q, A, \{0\}_N; 0) \\[0.5ex]
			& = I(C; B, \{0\}^m, A; 0) \,,
	\end{align*}
	as \( \sum_{b+p=N} (-1)^b \binom{N}{p} = (1 - 1)^N = \delta_{N = 0} \), and only the \( N = 0 \) term survives.  Then in terms of multiple polylogarithms (via \autoref{eqn:mpl}), after applying the affine invariance \autoref{eqn:affine}) this is just
	\[
		= \Li_{1,m+1}\Big(\frac{C}{B}, \frac{B}{A}\Big) \,.
	\]
	In particular, we will find
	\begin{align*}
		& \sum_{\substack{a+b=m \\ a,b \geq 0}} (-1)^b \cdot -\ii \big[ \text{\,summation in Eqn.~\eqref{eqn:Iab:part2} } \big]  \\
		& = -\ii \big[ -\Li_{1,m+1}(-\ii x, \ii) + \Li_{1,m+1}(-\ii x, -\ii) + \Li_{1,m+1}(\ii x, -\ii) - \Li_{1,m+1}(\ii x, \ii) \big] \,.
	\end{align*}
	(When computing \( [\bullet]_{x=1}^{x=-\ii} \), the bound \( x = -\ii \) should be taken as a limit; the apparent divergence here will cancel with another elsewhere later.)
	
	On the other hand, to deal with the terms in \autoref{eqn:Iab:part1}, we find (this is a shuffle-product identity, provable directly by induction) that,
	\begin{align*}
		I(C; B; \{0\}_a; 0) 
		&= \sum_{\substack{p+q = a \\ p,q \geq 0}} (-1)^p I(C; \{0\}_p, B; 0) I(0; \{0\}_q; C) \\[-0.5ex]
		&= \sum_{\substack{p+q = a \\ p,q \geq 0}} (-1)^{p+1} \Li_{p+1}\Big( \frac{C}{B} \Big) \cdot \frac{\log^q(C)}{q!} \,.
	\end{align*}
	In particular,
	\begin{align*}
		I(-\ii; 1, \{0\}_b; 0) - I(-\ii; -1, \{0\}_b; 0) &= \sum_{\substack{r+s=b \\ r,s \geq 0}} (-1)^{r+1} \Big( \Li_{r+1}(-\ii)  - \Li_{r+1}(\ii) \Big) \frac{\log^s(-\ii)}{s!}  \\[-0.5ex]
		& = 2 \ii  (-1)^b \sum_{\substack{r+s=b \\ r,s \geq 0}} \beta(r+1) \cdot \frac{1}{s!} \Big( \frac{\ii \pi}{2} \Big)^{s} \,,
	\end{align*}
	using that \( \Li_{r+1}(-\ii)  - \Li_{r+1}(\ii) = -2\ii \Im(\Li_{r+1}(\ii)) = -2\ii \beta(r+1) \), and  \( \log(-\ii) = -\frac{\ii \pi}{2} \).  Likewise
	\begin{align*}
		I(x; \ii, \{0\}_a; 0) - I(x; -\ii, \{0\}_a; 0)
		& = \sum_{\substack{p+q=a \\ p,q \geq 0}} (-1)^{p+1} \big( \Li_{p+1}(-\ii x) - \Li_{p+1}(\ii x)\big) \frac{\log^q(x)}{q!} \,.
	\end{align*}
	Then by gathering the terms with \( s + q = k \), and applying the binomial theorem, we find 
	\begin{align*}
	& \sum_{\substack{a+b=m \\ a,b\geq 0}} (-1)^b \cdot -\ii \big[ \text{$(I+I)\cdot(I+I)$ terms in Eqn.~\eqref{eqn:Iab:part1}} \big] \\[-1.0ex]
	& = \sum_{\substack{k+p+r=m \\ r,p,k \geq 0}} 2 \ii  (-1)^{1+p} \beta(r+1) \big( \Li_{p+1}(-\ii x) - \Li_{p+1}(\ii x)\big) \cdot \frac{1}{k!} \Big( \frac{\ii \pi}{2}  + \log(x) \Big)^{k} \,.
	\end{align*}
	
	\paragraph{Preliminary expression} From the above evaluations and simplifications, we find that
	\begin{align*}
		& \sum_{\substack{a+b=m \\ a,b \geq 0}} (-1)^b \widetilde{I_{a,b}} 
		 = 
		\begin{aligned}[t] 
		-\ii \bigg[ & \sum_{\substack{k+p+r=m \\ r,p,k \geq 0}} 2 \ii  (-1)^{1+p} \beta(r+1) \big( \Li_{p+1}(-\ii x) - \Li_{p+1}(\ii x)\big) \cdot \frac{1}{k!} \Big( \frac{\ii \pi}{2}  + \log(x) \Big)^{k} \\[-1ex]
		& -\Li_{1,m+1}(-\ii x, \ii) + \Li_{1,m+1}(-\ii x, -\ii) + \Li_{1,m+1}(\ii x, -\ii) - \Li_{1,m+1}(\ii x, \ii) \, \bigg] \,.
			\end{aligned} 
	\end{align*}
	At \( x = 1 \), this expression is finite, and simplifies to
	\begin{align*}
	 \sum_{\substack{a+b=m \\ a,b \geq 0}} (-1)^b \widetilde{I_{a,b}} \Big\rvert_{x=1} 
	& = 
	\begin{aligned}[t]
	-\ii \bigg[ 
	& \sum_{\substack{k+p+r=m \\ r,p,k \geq 0}} 4 (-1)^{p+1} \beta(r+1) \beta(p+1) \frac{1}{k!} \Big( \frac{\ii \pi}{2} \Big)^k \\[-1.5ex]
	& -\Li_{1,m+1}(-\ii, \ii) + \Li_{1,m+1}(-\ii, -\ii) + \Li_{1,m+1}(\ii, -\ii) - \Li_{1,m+1}(\ii, \ii) 
	 \bigg] \,.
	\end{aligned} 
	\end{align*}
	On the other hand, we need to look at \( x \to -\ii \) more carefully, because of the previously noted divergence, which arises from the terms \( \Li_1(\ii x) \), and from \( \Li_{1,m+1}(\ii x, \pm \ii) \), as \( x \to -\ii \).  By the stuffle product, we can write
	\[
		\Li_{1,m+1}(\ii x, \pm \ii) = \Li_{1}(\ii x) \Li_{m+1}(\pm \ii) - \Li_{m+1, 1}(\pm \ii, \ii x) - \Li_{m+2}(\mp x) \,,
	\]
	with the divergence as \( x \to -\ii \) confined now to \( \Li_1(\ii x) \), as well.  Rewriting the \( I_{a,b} \)-sum via this, and using \( \Li_{m+1}(-\ii) - \Li_{m+1}(\ii) = -2 \ii \beta(m+1) \), gives
		\begin{align}
	\nonumber & \hspace{-1.8em} \sum_{\substack{a+b=m \\ a,b \geq 0}} (-1)^b \widetilde{I_{a,b}} =  \\[0.5ex]
	&
	\label{eqn:Iab:asxtomi} \begin{aligned}[c] 
	\mathllap{-\ii \bigg[} & 
	 \sum_{\substack{k+p+r=m \\ r,p,k \geq 0}} 2 \ii  (-1)^{1+p} \beta(r+1) \big( \Li_{p+1}(-\ii x) - \delta_{p\geq1} \Li_{p+1}(\ii x)\big) \cdot \frac{1}{k!} \Big( \frac{\ii \pi}{2}  + \log(x) \Big)^{k} \\[-1.0ex]
	& + \Li_{1}(\ii x) \cdot \Big( \!\!\!	\sum_{\substack{r+k=m \\ r,k \geq 0}} 2 \ii \beta(r+1) \frac{1}{k!} \Big( \frac{\ii \pi}{2}  + \log(x) \Big)^{k}
	-2 \ii \beta(m+1) \Big) 
	 - \Li_{m+2}(x) + \Li_{m+2}(-x) \\[-0.5ex]
	& -\Li_{1,m+1}(-\ii x, \ii) + \Li_{1,m+1}(-\ii x, -\ii) 
 - \Li_{m+1, 1}(-\ii, \ii x)+ \Li_{m+1, 1}(\ii, \ii x) \, \smash{\bigg]} \,.
	\end{aligned} 
	\end{align}
	Since
	\[
		\lim_{x\to-\ii} \Li_1(\ii x) \cdot \Big( \frac{\ii\pi}{2} + \log(x) \Big) = 0 \,,
	\]
	and every term in the coefficient of \( \Li_1(\ii x) \) is divisible by \( \big( \frac{\ii\pi}{2} + \log(x) \big) \) (as the \( k = 0 \) term cancels against the additional \( -2\ii\beta(m+1) \)), we conclude the divergence in \( \sum_{a+b=m} (-1)^b I_{a,b} \) indeed cancels out as \( x \to -\ii \).  We can therefore formally set the term \( \Li_{1}(\ii x) = 0 \) in \autoref{eqn:Iab:asxtomi}, and substitute \( x = -\ii \) elsewhere to obtain the desired evaluation.  Namely (with abuse of notation via \( \delta_{p\geq1} \) for simplicity), and noting only the \( k = 0 \) term survives:
	\begin{align*}
		& \sum_{\substack{a+b=m \\ a,b \geq 0}} \!\! (-1)^b \widetilde{I_{a,b}} \rvert_{x=-\ii} = 
		 \begin{aligned}[t] 
		-\ii \bigg[ & 
		\sum_{\substack{r+p=m \\ r,p \geq 0}} \!\! 2 \ii  (-1)^{1+p} \beta(r+1) \big( \Li_{p+1}(-1) - \delta_{p\geq1} \Li_{p+1}(1)\big) 
		+ 2\ii \beta(m+2)  \\
		& -\Li_{1,m+1}(-1, \ii) + \Li_{1,m+1}(-1, -\ii) 
		- \Li_{m+1, 1}(-\ii, 1)+ \Li_{m+1, 1}(\ii, 1)  \smash{\bigg]} \,.
		\end{aligned} 
	\end{align*}

	\paragraph{Final expression} Assembling the above simplifications, and limiting behaviour as \( x \to -\ii \), we obtain the following expression for the multiple \( T \) value.  Notice that all the depth 2 terms appear in complex-conjugate pairs, which we have combined to write via their imaginary parts.
	\begin{Thm}\label{thm:Tm21m1:any}
		The following evaluation holds for all \( m \geq 1 \),
	\begin{align*}
	& T(\overline{2}, \{1\}_{m-1}, \overline{1}) 
		= \\
	&	\ii^{\delta_{\text{$m$\,even}}} \Big[ \, \begin{aligned}[t]
	& \!\! \sum_{\substack{k+p+r=m \\ r,p,k\geq 0}} \!\! 	4 \ii  (-1)^{p+1} \beta(r+1) \beta(p+1) \frac{1}{k!} \Big( \frac{\ii \pi}{2} \Big)^k  \\
	& +  \sum_{\substack{r+p=m \\ r,p \geq 0}} \!\! 2 (-1)^{1+p} \beta(r+1) \big( \Li_{p+1}(-1) - \delta_{p\geq1} \Li_{p+1}(1)\big) 
	+ 2\beta(m+2)  \\
	& - 2 \Im\bigl(
		\Li_{1,m+1}(-1,-\ii) + \Li_{m+1,1}(-\ii, 1) + \Li_{1,m+1}(\ii, -\ii) + \Li_{1,m+1}(-\ii,-\ii)
		\smash{\big) \Big]} \,.
		\end{aligned}
	\end{align*}
	\end{Thm}
	
	We are interested in the case \( m = 2\ell \) even.  In this case, we know that
	\[
		2 \ii \Im ( \Li_{1,2\ell+1}(-1,-\ii) ) = \Li_{1,2\ell+1}(-1,-\ii) - \Li_{1,2\ell+1}(-1,\ii) \,,
	\]
	must reduce to depth 1, since  the depth 2 parity theorem (e.g. \cite{panzer16Parity}, but also \cite[\S2.6]{goncharov01}) generally gives a reduction for the combination
	\[
		\Li_{a,b}(x,y) - (-1)^{a+b - 2}\Li_{a,b}(x^{-1}, y^{-1}) = \text{depth 1 \& products} \,,
	\]
	so our imaginary-part combination has the correct sign to reduce.  
	
	On the other hand, we can bypass such considerations entirely, as the MTV is necessarily real, so by taking the real-part of \autoref{thm:Tm21m1:any}, we obtain the following corollary directly.  (As \( \ii^{\delta_{\text{$m$\,even}}} = \ii \), we take the imaginary part of everything inside the square brackets, which also forces \( k \) to be even in the first line.)
	\begin{Cor}
		The following evaluation holds for all \( \ell \geq 1 \)
		\[
			T(\overline{2}, \{1\}_{2\ell-1}, \overline{1}) = 
			  \sum_{\substack{r+p+2k=2\ell \\ r,p,k \geq 0}} \!\! 4 (-1)^{r+k} \beta(r+1) \beta(p+1) \cdot \frac{1}{(2k)!} \Big( \frac{\pi}{2} \Big)^{2k} \,.
		\]
	\end{Cor}
	\noindent This concludes the proof of \autoref{thm:qn2:mtv}.\hfill\qedsymbol\medskip

	We note however, that in the case of \( m = 2\ell+1 \) odd, no particularly pleasant simplification of \autoref{thm:Tm21m1:any} seems possible.  One can show, for example,
	\[
		T(\overline{2},\overline{1}) = 8 \Im\Bigl(\Li_3\Bigl(\frac{1+\ii}{2}\Bigr)\Bigr) + 4 \beta(2) \log (2)-\frac{3 \pi ^3}{16}-\frac{\pi}{4}\log ^2(2) \,,
	\]
	which already begins to invoke new irreducible terms outside of Riemann zeta values and Dirichlet beta values, and even the alternating double zeta values.  Nevertheless, \autoref{thm:Tm21m1:any} always gives us an expression for \( T(\overline{2}, \{1\}_{2\ell}, \overline{1}) \), in terms of depth 2 level 4 coloured MZV's; we leave further investigation in this direction to the interested reader.
	
	\section{Proof of \autoref{thm:qn2:msv} (first part)} \label{sec:thm2b}
	
	The iterated integral representation of MSV's \cite[\S1.2]{xu2022alternating}, and \autoref{eqn:arctanint} above, allows us to write
	\begin{align*}
	S(\overline{2}, \{1\}_{m-1}, {1}) 
	&= (-1)^{\lfloor m/2 \rfloor} \int_0^1 \frac{\mathrm{d}t}{t} \circ \Big\{ \frac{-2 \,  \mathrm{d}t}{t^2 + 1} \Big\}_m \circ \frac{-2t\,\mathrm{d}t}{t^2+1}  \\
	&= \frac{1}{m!} (-1)^{\lfloor m/2 \rfloor} \int_{1 > s > r > 0}  (-2)^m (\arctan(s) - \arctan(r))^m \cdot \frac{\mathrm{d}s}{s}  \cdot \frac{-2r \, \mathrm{d}r}{r^2 + 1} 
	\end{align*}
 Make the same substitution from \autoref{eqn:int:subst} as before,
\[
(r,s) = \big( \tan(-\tfrac{1}{2\ii} \log(x)), \tan(-\tfrac{1}{2\ii} \log(y)) \big) \,,
\]
and we obtain
\[
S(\overline{2}, \{1\}_{m-1}, {1}) = \frac{1}{m!} (-1)^{\lfloor m/2 \rfloor} \cdot \ii^m \int_{\gamma \,,\, -\ii > y > x > 1} \frac{2(x-1)(\log(x) - \log(y))^m}{x(x+1)(y^2 - 1)} \mathrm{d}x \, \mathrm{d}y  \,,
\]
clockwise along the circular arc \( \gamma \) from \( 1 \) to \( -\ii \).  Note \( (-1)^{\lfloor m/2 \rfloor} \cdot \ii^m = \ii^{\delta_{\text{$m$\,odd}}}\) in this case.  Now we can utilise the partial fractions decomposition
\[
	\frac{x-1}{x(1+x)} = \frac{2}{x+1} - \frac{1}{x} \,,
\]
in order to compute the primitive with respect to \( x \), and proceed as before.  We find
\begin{align*}
	J_{a,b} &\coloneqq \int_{\gamma, -\ii>y>x>1} \frac{2}{a! \, b!} \frac{(x-1) \log(x)^a \log(y)^b}{x(1+x)(y^2 - 1)} \, \mathrm{d}x \, \mathrm{d}y \\[1ex]
	& = 
	 \bigg[ 
	\begin{aligned}[t] 
	& \,\, \big( 2 \, I(x; -1, \{0\}_a, 0) - I(x; 0; \{0\}_a; 0)\big) 
	\cdot \big( I(-\ii; 1, \{0\}_b; 0) - I(-\ii; -1, \{0\}_b; 0) \big) 
	\end{aligned} 
	\\[1.5ex]
	&   \begin{aligned}[c]
	\smash{\phantom{ = -\ii \bigg[ } +\sum_{\substack{p+q = a \\ p,q \geq 0}} \binom{b+p}{p}} \big( 
	& {-} 2 I(x; -1, \{0\}_q, 1, \{0\}_{b+p}; 0)
	+ 2 I(x; -1, \{0\}_q, -1, \{0\}_{b+p}; 0) \\[-1.5ex]
	& \quad \quad + I(x; 0, \{0\}_q, 1, \{0\}_{b+p}; 0)
	- I(x; 0, \{0\}_q, -1, \{0\}_{b+p}; 0)
	\big) \,\, {\bigg]_{x=1 \,.}^{x=-\ii}}  \end{aligned}
\end{align*}
Then using the same regularisation observation as in \autoref{sec:thm2a},  namely \autoref{eqn:Iab:sumreg}, and the equations thereafter, we can directly obtain the following.  (We note that all expressions are convergent this time.)
\begin{align*}
	&\sum_{\substack{a+b=m \\ a,b \geq 0}} (-1)^b J_{a,b} \, = \, 
	 \bigg[
	\begin{aligned}[t] 
		& \sum_{\substack{k+r=m \\ k, r \geq 0}} 2 \ii (-1)^1 \beta(r+1) \frac{1}{(k+1)!} \Big( \frac{\ii\pi}{2} + \log(x)\Big)^{k+1} \\
		& + \sum_{\substack{k+p+r = m \\ p,k,r \geq 0}} 2 \ii (-1)^{1+p} \beta(r+1) \cdot 2 \Li_{p+1}(-x) \cdot \frac{1}{k!} \Big( \frac{\ii \pi}{2} + \log(x) \Big)^k \\[-1.5ex]
		& + \Li_{m+2}(-x) - \Li_{m+2}(x) - 2 \Li_{1,m+1}(-x.-1) + 2 \Li_{1,m+1}(-x,1)
	{\,\bigg]_{x=1\,.}^{x=-\ii}} 
	\end{aligned}
\end{align*}
	Whence, after some straight-forward simplification (such as \( \big( \frac{\ii\pi}{2} + \log(-\ii)\big)^{k+1} = 0\), as \( k+1 > 0\)), we obtain the following.
	\begin{Thm}\label{thm:Sm2111:any}
	The following evaluation holds for all \( m \geq 1 \):
	\begin{align*}
		& S(\overline{2}, \{1\}_{m-1}, 1) = \\
		& \ii^{\delta_{\text{$m$\,odd}}} \bigg[
		\begin{aligned}[t]
			& \sum_{\substack{p+r=m \\ p,r \geq 0}} 2 \ii (-1)^{p+1} \beta(r+1) \cdot 2 \Li_{p+1}(\ii)
			+ \sum_{\substack{p+r=m \\ p,r \geq 0}} 2 \ii \beta(r+1) \cdot \frac{1}{(k+1)!} \Big( \frac{\ii \pi}{2} \Big)^{k+1} \\[-0.5ex]
			& - \sum_{\substack{k+p+r=m \\ p,k,r \geq 0}} 2 \ii (-1)^{p+1} \beta(r+1) \cdot 2 \Li_{p+1}(-1) \cdot \frac{1}{k!} \Big( \frac{\ii \pi}{2} \Big)^k  + 2 \ii \beta(m+2) - \Li_{m+2}(-1) \\[-1.5ex]
			& + \Li_{m+2}(1) - 2 \Li_{1,m+1}(\ii,-1) + 2 \Li_{1,m+1}(\ii,1) + 2 \Li_{1,m+1}(-1,-1) - 2 \Li_{1,m+1}(-1,1) \bigg] \,.
		\end{aligned}
	\end{align*}
	\end{Thm}
	Although it might be possible to still simplify this in arbitrary weight, for the purposes of \autoref{thm:qn2:msv} (first part), we are only interested in the case \( m = 2\ell \) even.  After taking the real part (for the MSV of interest is real), the depth 2 contribution is then
	\begin{equation}\label{eqn:Sv1:dp2part}
		\mathcal{D}_1\coloneqq 2 \Li_{1,2\ell+1}(-1,-1) - 2 \Li_{1,2\ell+1}(-1,1) + \Re\bigl( -2 \Li_{1,2\ell+1}(\ii,-1) + 2 \Li_{1,2\ell+1}(\ii, 1) \bigr) \,.
	\end{equation}
	To prove \autoref{thm:qn2:msv} (first part), it suffices to express \( R \) via \( \zeta(\overline{2\ell+1},1) \) and products of Riemann zeta values and Dirichlet beta values.  Let us do this term by term, after recalling a useful identity.
	
	\paragraph{The shuffle antipode} The following identity of iterated integrals is straight-forward to prove by induction (using the recursive definition of the shuffle product, c.f. \cite[Eqn.~(29)]{goncharovGalois05} or \cite[Lemma 4.2.1]{glanoisThesis16})
	\begin{equation}\label{eqn:shuffle:anti}
		\sum_{i=0}^n (-1)^i \int_a^b \omega_{1} \omega_{2} \ldots \omega_{i-1} \omega_i \cdot \int_a^b \omega_{n} \omega_{n-1} \cdots \omega_{i+2} \omega_{i+1} = 0 \,.
	\end{equation}
	By applying this to \( I(1; x_n, \ldots, x_1; 0) \), \( x_i \in \{ 0, \pm 1, \pm \ii \} \) we can obtain a number of relations involving level $1, 2$ or $4$ coloured MZV's (possibly under shuffle regularisation).

	\paragraph{Rewriting of $\Li_{1, 2\ell+1}(-1, 1)$ {\normalfont \&} $\Li_{1,2\ell+1}(\ii, 1)$} Applying the shuffle antipode gives 
	\[
	2 \Li_{1, 2\ell+1}(-1, 1) + \sum_{r=1}^{2\ell+1} (-1)^r \Li_r(-1) \Li_{2\ell+2-r}(-1) = 0 \,,
	\]
	as the integral representation  
	\[
		\Li_{1,2\ell+1}(-1, 1) = I(1; -1, \{0\}_{2\ell}, -1; 0)
	\]
	is palindromic in even weight, so the terms of maximal length combine.  Therefore \( \Li_{1,2\ell+1}(-1, 1) \) is already reducible to products of (alternating) Riemann zeta values.
	
	The same palindromicity holds for \[
		\Li_{1,2\ell+1}(\ii, 1) = I(1; -\ii, \{0\}_{2\ell}, -\ii; 0) \,,
	\]
	and we likewise obtain
	\[
		2 \Li_{1, 2\ell+1}(\ii, 1) + \sum_{r=1}^{2\ell+1} (-1)^r \Li_r(\ii) \Li_{2\ell+2-r}(\ii) = 0 \,,
	\]
	showing it is reducible to products of Riemann zeta values and Dirichlet beta values.
		
	\paragraph{Rewriting of $\Li_{1,2\ell+1}(-1,-1)$} Applying the shuffle antipode gives (under shuffle regularisation)
	\[
		\Li_{1,2\ell+1}(-1,-1) + \Li_{1,2\ell+1}(1, -1) + \sum_{r=1}^{2\ell+1} (-1)^r \Li_r(-1) \Li_{2\ell + 2 - r}(1) = 0 \,.
	\]
	For a single leading 1, the shuffle and stuffle regularisations are equal (c.f. \S13.3.1 and Theorem 13.3.9 in \cite{zhaoFunctions}, as well as \S13.3.2 for the proof), so we can directly compute
	\[
		\Li_{1,2\ell+1}(1, -1) = \Li_1(1) \Li_{2\ell+1}(-1) - \Li_{1,2\ell+1}(1, -1) - \Li_{2\ell+2}(-1) \,,
	\]
	and (with convergent expressions only) obtain
	\[
		\Li_{1,2\ell+1}(-1, -1) = \Li_{2\ell+1, 1}(-1, 1) + \Li_{2\ell+2}(-1) - \sum_{r=1}^{2\ell} (-1)^r \Li_r(-1) \Li_{2\ell+2-r}(1) \,.
	\]
	This is now in the required form.

	\paragraph{Rewriting of $\Li_{1,2\ell+1}(\ii, -1)$}	Finally, applying the shuffle antipode here gives
	\[
		\Li_{1, 2\ell+1}(\ii, -1) + \Li_{1, 2\ell+1}(-\ii, -1) + \sum_{r=1}^{2\ell+1} (-1)^r \Li_r(\ii) \Li_{2\ell+2-r}(-\ii) = 0 \,.
	\]
	In this case, we don't obtain a reduction directly, but notice
	\[
		\Li_{1, 2\ell+1}(\ii, -1) + \Li_{1, 2\ell+1}(-\ii, -1) = 2 \Re \big( \Li_{1, 2\ell+1}(\ii, -1)  \big) \,,
	\]
	and since the latter real part is contribution which actually appears in \( \mathcal{D}_1 \) (see \autoref{eqn:Sv1:dp2part} above) this is sufficient.
	
	\paragraph{Final result for \( S(\overline{2}, \{1\}_{2\ell-1}, 1) \)} With the above computations we have essentially proven \autoref{thm:qn2:msv} (first part).  It is however relatively straightforward to substitute the above evaluations into \autoref{thm:Sm2111:any}, simplify, and write the MPL's as MZV's.  After converting alternating single zeta values to \( \log(2) \) or Riemann zeta values, we obtain the following explicit result, which makes the form of the evaluation manifest.
	
	\begin{Cor}\label{cor:qn2:msvpart1}
		The following evaluation holds for any \( \ell \geq 1 \),
		\begin{align*}
			& S(\overline{2}, \{1\}_{2\ell-1},1)
			=  \\[1ex]
			& 2 \zeta(\overline{2\ell+1},1) + 2^{-1-2\ell} \zeta(2\ell+2) + \big(2^{1-2\ell} - 4 \big) \log(2) \zeta(2\ell+1) \\
			& + \!\! \sum_{\substack{r + 2k = 2\ell-1 \\ r,k \geq 0}} \!\! 4 \log(2)  \beta(r+1) \frac{(-1)^k}{(2k+1)!} \Big( \frac{\pi}{2} \Big)^{2k+1}
			 + \!\! \sum_{\substack{p+q=2\ell+2 \\ p, q \geq 2}} \!\! (-1)^p (1 - 2^{1-p}) (3 - 2^{1-q}) \zeta(p) \zeta(q) \\[-0.5ex]
			 & + \!\! \sum_{\substack{r + 2k = 2\ell \\ k, r \geq 0}} \!\! 2 \beta(r+1)  \, \frac{(-1)^{k+1}}{(2k+1)!} \Big( \frac{\pi}{2} \Big)^{2k+1} 
			  + \!\! \sum_{\substack{p+r= 2\ell \\ p, r \geq 0}} \!\! 2 (-1)^p \beta(r+1) \beta(p+1) \\[-0.5ex]
			  & + \!\! \sum_{\substack{p + r + 2k = 2\ell-1 \\ p \geq 1 ,\, r, k \geq 0}} \!\! 4 (-1)^{p+k} \beta(r+1) \, (1 - 2^{-p}) \zeta(p+1) \, \frac{1}{(2k+1)!} \Big( \frac{\pi}{2} \Big)^{2k+1} \,.
		\end{align*}
	\end{Cor}

	In particular, the form of \( S(\overline{2}, \{1\}_{2\ell-1}, 1) \) conjectured by Xu, Yan and Zhao \cite[Question 2, p.~21]{xu2022alternating} holds.  This completes the proof of \autoref{thm:qn2:msv} (first part). \hfill \qedsymbol
	
		\section{Proof of \autoref{thm:qn2:msv} (second part)} \label{sec:thm2c}

	Now, the iterated integral representation of MSV's \cite[\S1.2]{xu2022alternating}, plus \autoref{eqn:arctanint} and the substitution \autoref{eqn:int:subst} from before, allows us to write
	\begin{align*}
		S(\overline{2}, \{1\}^{m-1}, \overline{1}) &= (-1)^{\lfloor m/2 \rfloor} \int_0^1 \frac{\mathrm{d}t}{t} \circ \Big\{ \frac{-2 \, \mathrm{d}t}{t^2 + 1} \Big\}_m \circ \frac{-2 t \, \mathrm{d}t}{t^2 - 1} \\
		& = \frac{1}{m!} (-1)^{\lfloor m/2 \rfloor} \int_{1 > s > r > 0} (-2)^m (\arctan(s) - \arctan(r))^m \cdot \frac{\mathrm{d}s}{s}  \cdot \frac{-2r\,\mathrm{d}r}{r^2 - 1} \\ 
		& = \frac{\ii^{\text{$m$\,odd}}}{m!} \int_{\gamma,\,1 >  y > x > -\ii} \frac{4 (1-x) (\log(x) - \log(y))^m}{(x+1)(x^2 + 1) (y^2 - 1)}
	\end{align*}
	By the partial fraction expansion
	\[
		\frac{2(1-x)}{(1+x)(1 + x^2)} = \frac{2}{x - (-1)} - \frac{1}{x - \ii} - \frac{1}{x - (-\ii)} \,,
	\]
	and the regularisation simplification from \autoref{eqn:Iab:sumreg}, we likewise readily obtain
	\begin{align*}
		&S(\overline{2}, \{1\}_{m-1}, \overline{1}) = \\
		& \ii^{\text{$m$\,odd}} \Big[ \begin{aligned}[t]
		& \sum_{\substack{k+ p  + r = m \\ p,k,r\geq0}} 2 \ii (-1)^{p+1} \beta(r+1) \big( 2 \Li_{p+1}(-x) - \Li_{p+1}(-\ii x) - \Li_{p+1}(\ii x) \big) \frac{1}{k!} \Big( \frac{\ii\pi}{2} + \log(x) \Big)^k \\
		& - 2 \Li_{1,m+1}(-x, -1) + 2 \Li_{1, m+1}(-x, 1) + \Li_{1,m+1}(\ii x, -\ii) \\
		&  - \Li_{1,m+1}(\ii x, \ii) + \Li_{1,m+1}(-\ii x, \ii) - \Li_{1, m+1}(-\ii x, -\ii) 
		\Big]_{x=1 \,,}^{x=-\ii} \end{aligned}
	\end{align*}
	where, again, the lower bound \( x = -\ii \) must be interpreted via the limit; the apparent divergence cancels as before (but we leave the technicalities to the reader).  After rewriting to remove the apparent singularity (via the stuffle-product), and some straight-forward simplification (including using
	\[
		\Li_{p+1}(-\ii) + \Li_{p+1}(\ii) = {2^{-p}} \Li_{p+1}(-1) \,,
	\]
	the distribution relation of level 2), we obtain the following result.
	
	\begin{Thm}\label{thm:Sm2111m1:any}
	The following evaluation holds for all \( m \geq 1 \),
	\begin{align}
		\nonumber & \hspace{-2em} S(\overline{2}, \{1\}_{m-1},\overline{1}) = \\
		 \ii^{\text{$m$\,odd}} \Big[ 
		\nonumber & 2 \ii \beta(m+2) + \sum_{\substack{k+p+r = m \\ k,p,r \geq 0}} 2 \ii (-1)^p \beta(r+1) \cdot  \big(2 - 2^{-p} \big) \Li_{p+1}(-1) \frac{1}{k!} \Big( \frac{\ii\pi}{2} \Big)^{k} \\
		\nonumber & + \sum_{\substack{p+r = m \\ p,r \geq 0}} 2 \ii (-1)^{p+1} \beta(r+1) \cdot (2 \Li_{p+1}(\ii) - \delta_{p\geq1} \Li_{p+1}(1) - \Li_{p+1}(-1) )  \\
		\nonumber & + 2 \Li_{1,m+1}(\ii, 1) - 2 \Li_{1,m+1}(-1, 1) + 2 \Li_{1, m+1}(-1, -1)  - 2 \Li_{1, m+1}(\ii, -1) \\ 
		\label{eqn:penultimate} &+ \big( \Li_{1,m+1}(-\ii, -\ii) + \Li_{1, m+1}(-1, \ii) \big) + \bigl( -\Li_{1,1+m}(\ii, -\ii) + \Li_{m+1,1}(\ii, 1) \big) \\
		\nonumber & + \bigl( -\Li_{m+1,1}(-\ii, 1) - \Li_{1, m+1}(-1, -\ii) - \Li_{1, m+1}(-\ii, \ii) + \Li_{1,m+1}(\ii, \ii) \big) \Big] \,.
	\end{align}
	\end{Thm}	

	As before, we are only interested in the case \( m = 2\ell \) even, and can apply the real-part to this evaluation in order to obtain some quick simplifications.  In particular, the penultimate line \eqref{eqn:penultimate} satisfies
	\begin{align*}
	& \Re\big[ \big( \Li_{1,m+1}(-\ii, -\ii) + \Li_{1, m+1}(-1, \ii) \big) + \big( -\Li_{1,1+m}(\ii, -\ii) + \Li_{m+1,1}(\ii, 1) \big) \big] \\
	& = \Re\big[ \big( \Li_{1,m+1}(\ii, \ii) + \Li_{1, m+1}(-1, -\ii) \big) + \big( -\Li_{1,1+m}(-\ii, \ii) + \Li_{m+1,1}(-\ii, 1) \big) \big] \,;
	\end{align*}
	from this we can cancel/combine many of these terms with those in (the real part of) the final line.  We obtain
	\begin{align*}
	& S(\overline{2}, \{1\}_{2\ell-1},\overline{1}) = \\
		& \Big[ \begin{aligned}[t]
		& \sum_{\substack{p+r+2k = 2\ell-1 \\ k,p,r\geq0}} 2 (-1)^{1+k+p} \beta(r+1) \cdot  \big(2 - 2^{-p} \big) \Li_{p+1}(-1) \frac{1}{(2k+1)!} \Big( \frac{\pi}{2} \Big)^{2k+1} \\
		& + \sum_{\substack{p+r = 2\ell \\ p,r \geq 0}} 4 (-1)^{p} \beta(r+1) \beta(p+1) - 2 \Li_{1,2\ell+1}(-1, 1) + 2 \Li_{1, 2\ell+1}(-1, -1)   \\
		& + \Re\big( 2 \Li_{1,2\ell+1}(\ii, 1) - 2 \Li_{1, 2\ell+1}(\ii, -1) \big) +  \Re\bigl( -2\Li_{1,2\ell+1}(-\ii, \ii) + 2\Li_{1, 2\ell+1}(\ii,\ii) \big) \Big] \,.
		\end{aligned}
		\end{align*}
		From the shuffle-antipode relations (see \autoref{eqn:shuffle:anti}, and the deductions thereafter, plus some regularisation steps as necessary), we have (again)
		\begin{align*}
			\Li_{1,2\ell+1}(-1,-1) &= \begin{aligned}[t] 
			& \Li_{2\ell+1,1}(-1,1) + \Li_{2\ell+2}(-1)  - \sum_{r=1}^{2\ell} (-1)^r \Li_r(-1) \Li_{2\ell+2-r}(1) \,, \end{aligned} \\
			2\Li_{1,2\ell+1}(-1,1) &= - \sum_{r=1}^{2\ell+1} (-1)^r \Li_r(-1) \Li_{2\ell+1-r}(-1) \,, \\
			2 \Li_{1,2\ell+1}(\ii, 1) &= - \sum_{r=1}^{2\ell+1} (-1)^i \Li_r(\ii) \Li_{2\ell+2-r}(\ii) \\[-1ex]
			\Re \big( \Li_{1,2\ell+1}(\ii,-1) \big) & = \Li_{1,2\ell+1}(\ii, -1) + \Li_{1,2\ell+1}(-\ii, -1) = - \sum_{r=1}^{2\ell+1} (-1)^r \Li_r(\ii) \Li_{2\ell+2-r}(-\ii)			
		\end{align*}
		Hence all of the depth 2 terms are reduced to products and a single \( \zeta(\overline{2\ell+1},1) \), except for the final combination
		\[
			\mathcal{D}_2 \coloneqq \Re( -2 \Li_{1,2\ell+1}(-\ii, \ii) + 2 \Li_{1,2\ell+1}(\ii, \ii)) \,.
		\]
		We can now make the following observation (directly by their definitions), that
		\[
			\mathcal{D}_2 = - M_{-1,-1}^{\od,\od}(1, 2\ell+1) = -4 t(\overline{1},\overline{2\ell+1})  = -4 \sum_{m_1 > m_2 > 0} \frac{(-1)^{m_1} (-1)^{m_2}}{(2m_1-1) (2m_2-1)^{2\ell+1}}\,,
		\]
		so we can express \( S(\overline{2}, \{1\}_{2\ell-1}, \overline{1}) \) in terms of an alternating version of some depth 2 multiple \( t \) value.  After simplification, and rewriting MPL's as (alternating MZV's), we obtain the following.
		
		\begin{Thm}\label{thm:Sm2111m1:even}
			The following evaluation holds for all \( \ell \geq 1 \),
			\begin{align*}
			& S(\overline{2}, \{1\}_{2\ell-1}, \overline{1}) = \\
			& \sum_{\substack{p+r+2k = 2\ell-1 \\ k,p,r\geq 0}} 2 (-1)^{1+k+p} \beta(r+1) \cdot  \big(2 - 2^{-p} \big) \zeta(\overline{p+1}) \frac{1}{(2k+1)!} \Big( \frac{\pi}{2} \Big)^{2k+1} \\[-0.5ex]
			& + \sum_{\substack{p+r = 2\ell \\ p,r \geq 0}} 2 (-1)^{p} \beta(r+1) \beta(p+1)  
			+ \sum_{\substack{p+q = 2\ell+2 \\ p, q \geq 1}}  (-1)^p \zeta(\overline{p}) \zeta(\overline{q})  \\
			& + \sum_{\substack{p+q = 2\ell+2 \\ p \geq 1, q \geq 2}}   2(-1)^{p+1} \zeta(\overline{p}) \zeta(q) 
			\,\, + \,\, 2 \zeta(\overline{2\ell+2})  + 2 \zeta(\overline{2\ell+1}, 1) - 4 t(\overline{1},\overline{2\ell+1}) \,.
			\end{align*}
		\end{Thm}
	
		Our task now is to reduce \( t(\overline{1}, \overline{2\ell+1}) \) to alternating MZV's.
		
		\paragraph{Evaluation of \( t(\overline{1},\overline{2\ell+1}) \)}  The key input to this evaluation is a variation of the so-called generalised doubling identity.  (Note: a different application of generalised doubling was used in \cite{CK2242} to establish a reduction of \( \zeta(\overline{\ev}, \overline{\ev}) \) to non-alternating double zeta values; it would be interesting to investigate whether \( t(\overline{\od}, \overline{\od}) \) can similarly be reduced to alternating double zeta values, or non-alternating multiple $t$ values.)  By following the proof from \cite[\S14.2.5]{zhaoFunctions} or \cite[\S4]{mzvDM}, one can easily adapt it to obtain the following functional version.
		\begin{Lem}[`Doubly' generalised doubling relation]
			The following identity hold
			\begin{align*}
				& \frac{1}{2} \big( \Li_{s,t}(x,y) + \Li_{s,t}(-x,-y) \big) = \\
				& \sum_{j=1}^t \frac{1}{2^{s+t-j}} A_j^{t,s} \Li_{s+t-j,j}\bigl(x^2, \tfrac{y}{x} \bigr) + \sum_{j=1}^s  \frac{1}{2^{s+t-j}} A_i^{s,t} \Li_{i,s+t-j}\bigl(\tfrac{x}{y}, y^2\bigr) \\[-0.5ex]
				& - \sum_{j=1}^s \frac{1}{2} A_j^{s,t} \big( \Li_{j,s+t-j}\bigl( \tfrac{x}{y}, -y \bigr) + \Li_{j,s+t-j}\bigr( \tfrac{x}{y}, y \bigr) \big) - \frac{(s+t+1)!}{(s-1)! \, t!} 2^{-s-t} \Li_{s+t}(x^2) \,,
			\end{align*}
			where \( A_{j}^{s,t} = \binom{s+t-j-1}{t-1} \)\,.
		\end{Lem}
	
		If we specialise to the case \( s = 1, t = m \), \( x = y = \ii \), (and apply the stuffle product, to see the divergences cancel), we obtain
		\begin{equation}\label{eqn:dup:ii}
		\begin{aligned}[c]
			& \Li_{1,m+1}(\ii, \ii) + \Li_{1,m+1}(-\ii, -\ii) - \Li_{m+1,1}(-\ii, 1) - \Li_{m+1,1}(\ii, 1) \\
			& + \frac{1}{2^m} \Li_{m+1,1}(-1,1) + \frac{1}{2^m} \Li_{m+2}(-1) - 2 \sum_{j=1}^{m+1} \frac{1}{2^j} \Li_{j,m+2-j}(-1, 1) = 0 \,,
		\end{aligned}
		\end{equation}
		where the first line can be rewritten as \( \Re( \Li_{1,m+1}(\ii, \ii) - \Li_{m+1}(-\ii, 1) ) \).  We then invoke the shuffle-antipode once more to obtain (with shuffle-regularisation)
		\[
			\Li_{2\ell+1,1}(1,-\ii) + \Li_{1,2\ell+1}(-\ii,\ii) + \sum_{r=1}^{2\ell+1} (-1)^r \Li_r(-\ii) \Li_{2\ell+2-r}(1) \,,
		\]
		whence, (after converting to stuffle-regularisation)
		\begin{equation}\label{eqn:mi1:antipode}
		\begin{aligned}[c]
			& \Li_{2\ell+1,1}(-\ii, 1) = {}
			 \Li_{1,2\ell+1}(-\ii,\ii)  -\Li_{2\ell+2}(-\ii)  
			 + \sum_{r=1}^{2\ell} (-1)^i \Li_r(-\ii) \Li_{2\ell+2-r}(1) \,.
		\end{aligned}
		\end{equation}
		Combining \autoref{eqn:dup:ii} and \autoref{eqn:mi1:antipode}, we obtain the following result
		
		\begin{Prop}\label{prop:tm1modd:eval}
			The following evaluation holds for all \( \ell \geq 1 \),
		\begin{align*}
			 4 t(\overline{1},\overline{2\ell+1}) = {} 
			 & \frac{1}{2^{2\ell}} \zeta(\overline{2\ell+1},1) - 2 \sum_{j=1}^{2\ell+1} \frac{1}{2^j} \zeta(\overline{j}, 2\ell+2-j)   \\
			 & + \frac{3}{2^{2\ell+1}} \zeta(\overline{2\ell+2})
			- \sum_{r=1}^{2\ell} \frac{(-1)^r}{2^{r-1}} \zeta(\overline{r}) \zeta(2\ell+2-r) \,.
		\end{align*}
		In particular \( t(\overline{1}, \overline{2\ell+1}) \) is expressed through depth 2 alternating MZV's and products involving Riemann zeta values and \( \log(2) \).
		\end{Prop}
		
		\paragraph{Final result for \( S(\overline{2}, \{1\}_{2\ell-1}, \overline{1}) \)}  Combining \autoref{prop:tm1modd:eval} with \autoref{thm:Sm2111m1:even} readily leads to the following evaluation.  Here we have converted the alternating single zeta values to Riemann zeta values, or \( \log(2) \) in order to make the form of the evaluation for \autoref{thm:qn2:msv} (second part) clear.
		
		\begin{Cor}\label{cor:qn2:msvpart2}
		The following evaluation holds for any \( \ell \geq 1 \),
		\begin{align*}
			& S(\overline{2}, \{1\}_{2\ell-1}, \overline{1}) =  \\[1ex]
			& \sum_{\substack{ p + r + 2k = 2\ell-1 \\ k, r \geq 0 \,, p \geq 1}} 2 (-1)^{k+p} \beta(r+1) \, (1 - 2^{-p})(2 - 2^{-p}) \zeta(p+1) \, \frac{1}{(2k+1)!} \Big( \frac{\pi}{2} \Big)^{2k+1} \\
			& + \sum_{\substack{p+r = 2\ell \\ p,r \geq 0}}  \!\! 2(-1)^p \beta(r+1)\beta(p+1) + \sum_{\substack{p+q = 2\ell+2 \\ p, q \geq 2}} (-1)^{p} (1 - 2^{1-p})(3 - 2^{1-p} - 2^{1-q}) \zeta(p)\zeta(q) \\[0.5ex]
			& + \sum_{\substack{r+2k = 2\ell-1 \\ r, k \geq 0}} 2 \log(2) \beta(r+1) \, \frac{(-1)^k}{(2k+1)!} \Big( \frac{\pi}{2} \Big)^{2k+1} 
			 + (3\cdot2^{-2\ell-1} - 2) (1 - 2^{-2\ell-1}) \zeta(2\ell+2) \\
			& + (2^{1-2\ell} - 3) \log(2) \zeta(2\ell+1)  + (2 - 2^{-2\ell}) \zeta(\overline{2\ell+1},1) + 2  \cdot  \sum_{\substack{p + q = 2\ell + 2 \\ p, q \geq 1}} \frac{1}{2^p} \zeta(\overline{p}, q) \,.
		\end{align*}
		\end{Cor}
	
		\noindent This completes the proof of \autoref{thm:qn2:msv} (second part). \hfill \qedsymbol\medskip
	
		\paragraph{Concluding remarks} We point out, in the case \( \ell = 3 \), the following evaluation of the weighted sum
		\begin{align*}
			& \sum_{\substack{p + q = 8 \\ p, q \geq 1}} \frac{1}{2^p} \zeta(\overline{p}, q)  = 
			 -\frac{177}{1216} \zeta(\overline{7},1) -\frac{17}{304}\zeta(\overline{5},3) - \frac{275763}{3404800}  \zeta(2)^4 +\frac{1305}{4864} \zeta(3) \zeta(5) +\frac{63}{64}  \zeta(7) \log (2)\,,
		\end{align*}
		where the latter is obtained using the alternating MZV Data Mine \cite{mzvDM}.  In particular, since \( \zeta(\overline{5},3) \) is a necessary new algebra generator in weight 8 (assuming the standard MZV conjectures), we already cannot express this result, and thence neither \( t(\overline{1}, \overline{7}) \) nor \( S(\overline{2}, \{1\}_5, \overline{1}) \), via \( \zeta(\overline{7}, 1) \) alone.  Likewise, additional alternating MZV's of the form \( \zeta(\overline{2a+1}, 2b+1) \), \( a > b \) enter the picture as \( \ell \) increases.  This answers the final part of Question 2 \cite[p.~21]{xu2022alternating} \emph{negatively} (modulo the standard MZV conjectures). \medskip
		 
		For the interested reader, here are the full evaluations (in the MZV Data Mine basis), for \( \ell = 3, 4 \).  
		For  \( \ell = 5 \), weight 12, one would need to invoke three alternating double zeta values \( \zeta(\overline{11}, 1), \zeta(\overline{9},3) \) and \( \zeta(\overline{7},5) \) to give a similar evaluation. {
		\begin{align*}
			S(\overline{2}, \{1\}^5, \overline{1}) \,=\, {} 
			&
			\frac{559}{304} \zeta(\overline{7},1)
			-\frac{17}{304} \zeta (\overline{5}, 3)
			-4 \beta (2) \beta (6)
			-2 \beta (4)^2
			+\frac{465}{256} \pi  \beta (2) \zeta (5)
			\\
			&
			+\frac{21}{16} \pi  \beta (4) \zeta (3)
			-\frac{7}{128} \pi ^3 \beta (2) \zeta(3)
			+\pi  \beta (6) \log (2)
			-\frac{1}{24} \pi ^3 \beta (4) \log (2)
			\\
			&
			+\frac{1}{1920}\pi ^5 \beta	(2) \log (2)
			-\frac{10377}{2432}\zeta (3) \zeta (5)
			-\frac{127}{64} \zeta (7) \log (2)
			+\frac{84869 \pi^8}{326860800} \,, \\[3ex]
			S(\overline{2}, \{1\}^7, \overline{1}) \,=\, {} 
			& 
			\frac{47483}{25328}  \zeta (\overline{9},1)
			-\frac{3165}{101312}  \zeta (\overline{7},3)
			+\frac{8001 }{4096} \pi  \beta (2) \zeta (7)
			-\frac{155}{2048}  \pi ^3 \beta (2) \zeta (5)
			\\
			& 
			+\frac{465}{256} \pi  \beta (4) \zeta (5)
			+\frac{21}{16} \pi  \beta (6) \zeta(3)
			-\frac{7}{128} \pi ^3 \beta (4) \zeta (3)
			+\frac{7 }{10240} \pi ^5 \beta (2) \zeta (3)
			\\
			& 
			-4 \beta (2) \beta (8)
			-4 \beta (4) \beta (6)
			-\frac{1}{24} \pi ^3 \beta (6) \log (2)
			+\frac{1}{1920} \pi ^5 \beta (4) \log (2)
			\\
			& 
			-\frac{1}{322560} \pi ^7 \beta (2) \log(2)
			+\pi  \beta (8) \log (2)
			-\frac{511}{256} \zeta (9) \log (2)
			-\frac{3606645}{810496}  \zeta (3) \zeta (7)
			\\
			& 
			-\frac{33075465}{12967936}  \zeta(5)^2
			+\frac{1364516407 \pi ^{10}}{38822888079360} \,.%
		\end{align*}%
	}%
	\bibliographystyle{habbrv2}%
	\bibliography{STpm2111pm1}

\end{document}